\documentclass{amsart}
\usepackage[utf8]{inputenc}
\usepackage{a4}
\usepackage{graphics}
\usepackage{color}
\usepackage{amssymb}
\usepackage[all]{xy}
\usepackage{amsmath}
\usepackage{stmaryrd}
\usepackage{subfigure}
\usepackage{longtable}
\usepackage{setspace}
\usepackage{booktabs}

\usepackage{pst-plot}
\usepackage{pstricks}
\usepackage{color}

\CompileMatrices

\newtheorem{theorem}{Theorem}[section]
\newtheorem{lemma}[theorem]{Lemma}
\newtheorem{proposition}[theorem]{Proposition}
\newtheorem{corollary}[theorem]{Corollary}
\theoremstyle{definition}
\newtheorem{definition}[theorem]{Definition}
\newtheorem{example}[theorem]{Example}

\newtheorem{remark}[theorem]{Remark}
\theoremstyle{remark}

\numberwithin{equation}{section}
\def\Chi{{\mathbb X}}

\def\div{{\rm div}}

\def\quot{/\!\!/}
\def\mal{\! \cdot \!}

\def\rq#1{\widehat{#1}}
\def\t#1{\widetilde{#1}}
\def\b#1{\overline{#1}}
\def\bangle#1{\langle #1 \rangle}

\def\KK{{\mathbb K}}
\def\TT{{\mathbb T}}
\def\ZZ{{\mathbb Z}}

\def\QQ{{\mathbb Q}}
\def\PP{{\mathbb P}}
\def\Of{{\mathcal{O}}}

\def\WDiv{\operatorname{WDiv}}

\def\ord{{\rm ord}}

\def\Cl{\operatorname{Cl}}

\def\Pic{\operatorname{Pic}}

\def\Hom{{\rm Hom}}

\def\Spec{{\rm Spec}}
\def\Proj{{\rm Proj}}

\def\Pol{{\rm Pol}}

\def\cone{{\rm cone}}
\def\lin{{\rm lin}}

\DeclareMathOperator{\tail}{tail}

\DeclareMathOperator{\syz}{Rel}

\newcommand{\D}{\mathcal{D}}

\newcommand{\fan}{\Xi}

\newcommand{\xrays}{\fan^\times}
\newcommand{\xvers}{\fan^\times}

\newgray{lgray}{0.68}
\newgray{grayi}{0.68}
\newgray{grayii}{0.60}
\newgray{grayiii}{0.55}
\newgray{grayiv}{0.50}
\newgray{grayv}{0.45}
\newgray{grayvi}{0.40}

\newcommand{\zline}[1]{%
  \psline[linewidth=0.4pt,linestyle=dotted,dotsep=1pt](#1,-0.2)(#1,-.5)
}

\newcommand{\zlab}[1]{%
  \psline[linewidth=0.4pt,linestyle=dotted,dotsep=1pt](#1,-0.2)(#1,-.7)
  \fontsize{5}{6}%
  \rput(#1,-1.1){$#1$}
}

\newcommand{\mygrid}{
\zlab{-3}
\zlab{-2}
\zlab{-1}
\zlab{0}
\zlab{3}
\zlab{2}
\zlab{1}
\zline{-2.5}
\zline{-1.5}
\zline{-0.5}
\zline{0.5}
\zline{2.5}
\zline{1.5}
}

\newcounter{itemnumber}

\begin{document}
\title[The Cox ring of a $T$-variety]%
{The Cox ring of an algebraic \\
variety with torus action}
\author[J.~Hausen]{J\"urgen Hausen} 
\address{Mathematisches Institut, Universit\"at T\"ubingen,
Auf der Morgenstelle 10, 72076 T\"ubingen, Germany}
\email{hausen@mail.mathematik.uni-tuebingen.de}
\author[H.~S\"uss]{Hendrik S\"uss }
\address{Institut f\"ur Mathematik,
        LS Algebra und Geometrie,
        Brandenburgische Technische Universit\"at Cottbus,
        PF 10 13 44, 
        03013 Cottbus, Germany}
\email{suess@math.tu-cottbus.de}

\begin{abstract}
We investigate the Cox ring of a normal complete variety $X$
with algebraic torus action.
Our first results relate the Cox ring of $X$ to that
of a maximal geometric quotient of $X$.
As a consequence, we obtain a complete description
of the Cox ring in terms of generators and relations
for varieties with torus action of complexity one.
Moreover, we provide a combinatorial approach to the
Cox ring using the language of polyhedral divisors. 
Applied to smooth $\KK^*$-surfaces, our results
give a description of the Cox ring in terms 
of Orlik-Wagreich graphs.
As examples, we explicitly compute the  
Cox rings of all Gorenstein del Pezzo $\KK^*$-surfaces 
with Picard number at most two
and the Cox rings of projectivizations 
of rank two vector bundles as well as cotangent bundles
over toric varieties in terms of 
Klyachko's description.
\end{abstract}

\maketitle

\section{Introduction}

Let $X$ be a normal complete algebraic variety 
defined over some algebraically closed 
field~$\KK$ of characteristic zero
and suppose 
that the divisor class group $\Cl(X)$
is finitely generated.
The Cox ring of~$X$ is the graded 
$\KK$-algebra 
\begin{eqnarray*}
\mathcal{R}(X)
& = & 
\bigoplus_{\Cl(X)} \Gamma(X, \mathcal{O}_X(D)),
\end{eqnarray*}
see Section~\ref{coxring} for a detailed 
reminder.
A basic problem is to present $\mathcal{R}(X)$
in terms of generators and relations.
Besides the applications in number theory, 
see e.g.~\cite{DerTsch},
the knowledge of  generators and relations
also opens a combinatorial 
approach to geometric properties of 
$X$, see~\cite{BeHa2} and~\cite{Ha2}.

In the present paper, we investigate the 
case that $X$ comes
with an effective algebraic torus 
action $T \times X \to X$.
Our first result relates the 
Cox ring of $X$ to that of a 
maximal orbit space of the $T$-action.
For a point $x \in X$, 
denote by $T_x \subseteq T$ 
its isotropy group
and consider the non-empty 
$T$-invariant open subset
$$ 
X_0
\ := \ 
\{x \in X; \; \dim(T_x) = 0\}
\ \subseteq \ 
X.
$$
There is a geometric quotient 
$q \colon X_0 \to X_0 / T$ 
with an irreducible normal but possibly 
non-separated orbit space 
$X_0 / T$, see~\cite{Su},
and also for $X_0/T$ one can define 
a Cox ring.
Denote by $E_1, \ldots, E_m \subseteq X$
the ($T$-invariant) prime divisors 
supported in $X \setminus X_0$ and
by $D_1, \ldots, D_n \subseteq X$ 
those $T$-invariant prime divisors
who have a finite generic 
isotropy group of order $l_j > 1$.
Moreover, let $1_{E_k}$ and $1_{D_j}$ 
denote the canonical sections of the
divisors $E_k$ and $D_j$ respectively, 
and let $1_{q(D_j)} \in \mathcal{R}(X_0/T)$ 
be the canonical section of $q(D_j)$.

\begin{theorem}
\label{fingenchar}
There is a graded injection 
$q^* \colon \mathcal{R}(X_0/T) \to \mathcal{R}(X)$
of Cox rings and the assignments $S_k \mapsto 1_{E_k}$
and $T_j \mapsto 1_{D_j}$
induce an isomorphism of $\Cl(X)$-graded rings
\begin{eqnarray*}
\mathcal{R}(X)
& \cong &
\mathcal{R}(X_0/T)[S_1, \ldots, S_m, T_1, \ldots, T_n] 
\ / \ \bangle{T_j^{l_j} - 1_{q(D_j)}; \; 1 \le j \le n},
\end{eqnarray*}
where $\Cl(X)$-grading on the right hand side 
is defined by associating to $S_k$ the class of $E_k$
and to $T_j$ the class of $D_j$.
In particular, $\mathcal{R}(X)$ is finitely
generated if and only if 
$\mathcal{R}(X_0/T)$ is so. 
\end{theorem}

\goodbreak

If the dimension of $T$ equals that of $X$,
then our $X$ is a toric variety,
the subset $X_0 \subseteq X$ is
the open $T$-orbit, 
the divisors $E_1, \ldots, E_m$ 
are the invariant prime divisors of $X$ and there
are no divisors $D_j$.
Thus, for toric varieties, the above Theorem 
shows that the Cox ring is the polynomial
ring in the canonical sections of 
the invariant prime divisors and hence 
gives the result obtained by D.~Cox in~\cite{Co}.

The Cox ring $\mathcal{R}(X)$ can be further
evaluated by using the fact that $X_0/T$
admits a separation, i.e., a rational map 
$\pi \colon X_0/T \dashrightarrow Y$ to a variety $Y$,
which is a local isomorphism 
in codimension one.
After suitably shrinking, we may assume 
that there are prime divisors 
$C_0, \ldots, C_r$ on $Y$ 
such that each inverse image $\pi^{-1}(C_i)$
is a disjoint union of prime divisors $C_{ij}$,
where $1 \le j \le n_i$, the map
$\pi$ is an isomorphism over 
$Y \setminus (C_0 \cup \ldots \cup C_r)$
and all the $D_j$ occur among the $D_{ij} := q^{-1}(C_{ij})$.
Let $l_{ij} \in \ZZ_{\ge 1}$ denote the order
of the generic isotropy group of $D_{ij}$.

\begin{theorem}
\label{fingenchar2}
There is a graded injection 
$\mathcal{R}(Y) \to \mathcal{R}(X)$ 
of Cox rings
and the assignments $S_k \mapsto 1_{E_k}$
and $T_{ij} \mapsto 1_{q^{-1}(D_{ij})}$
induce an isomorphism of $\Cl(X)$-graded rings
\begin{eqnarray*}
\mathcal{R}(X)
& \cong &
\mathcal{R}(Y)[S_1, \ldots, S_m, T_{ij}; \; 0 \le i \le r, \, 1 \le j \le n_i]
\ / \ 
\bangle{T_i^{l_i} - 1_{C_i}; \; 0 \le i \le r}.
\end{eqnarray*}
where $T_i^{l_i} := T_{i1}^{l_{i1}} \cdots T_{in_i}^{l_{in_i}}$,
and the $\Cl(X)$-grading on the right hand side 
is defined by associating to $S_k$ the class of $E_k$ 
and to $T_{ij}$ the class of $D_{ij}$.
In particular, $\mathcal{R}(X)$ is finitely
generated if and only if 
$\mathcal{R}(Y)$ is so. 
\end{theorem}

Now suppose that  the $T$-action on $X$ 
is of complexity one, 
i.e., its biggest $T$-orbits are
of codimension one in $X$.
Then $X_0/T$ is of dimension one and 
has a separation $\pi \colon X_0/T \to \PP_1$.
Choose $r \ge 1$ and $a_0, \ldots, a_r \in \PP_1$ 
such that $\pi$ is an isomorphism over 
$\PP_1 \setminus \{a_0, \ldots, a_r\}$
and all the divisors $D_j$ occur among 
the $D_{ij} := q^{-1}(y_{ij})$, where 
$\pi^{-1}(a_i) = \{y_{i1}, \ldots, y_{in_i}\}$.
Let $l_{ij} \in \ZZ_{\ge 1}$ denote
the order of the generic isotropy group of 
$D_{ij}$.
For every $0 \le i \le r$,
define a monomial
$$
f_i 
\  := \
T_{i1}^{l_{i1}} \cdots T_{in_i}^{l_{in_i}}
\ \in \
\KK[T_{ij}; \; 0 \le i \le r, \; 1 \le j \le n_i].
$$
Moreover, write $a_i = [b_i,c_i]$ with 
$b_i,c_i \in \KK$ and 
for every $0 \le i \le r-2$ set 
$k = j+1 = i+2$ 
and define a trinomial
\begin{eqnarray*}
g_i 
& := &
(c_kb_j - c_jb_k)f_i
\ + \ 
(c_ib_k - c_kb_i)f_j
\ + \ 
(c_jb_i - c_ib_j)f_k.
\end{eqnarray*}

\begin{theorem}
\label{complexity1}
Let $X$ be a normal complete variety 
with 
finitely generated divisor class group
and an effective algebraic torus action 
$T \times X \to X$ 
of complexity one.
Then, in terms of the data defined
above, the Cox ring of $X$ is given as
\begin{eqnarray*}
\mathcal{R}(X)
&  \cong &
\KK[S_1,\ldots, S_m, T_{ij}; \; 0 \le i \le r, \; 1 \le j \le n_i] 
\ / \
\bangle{g_i; \; 0 \le i \le r-2},
\end{eqnarray*}
where $1_{E_k}$ corresponds to $S_k$,
and $1_{D_{ij}}$ to $T_{ij}$,
and the $\Cl(X)$-grading on the right hand side 
is defined by associating to $S_k$ the class 
of $E_k$ 
and to $T_{ij}$ the class of $D_{ij}$.
In particular, $\mathcal{R}(X)$ is finitely
generated. 
\end{theorem}

Note that finite generation of the Cox ring
for a complexity one torus action with $X_0/T$ 
rational may as well be deduced from~\cite{Kn}.

In Section~\ref{sec:crviapdiv}, we combine the 
results just presented with the description of 
torus actions in terms of polyhedral divisors
given in~\cite{MR2207875} and~\cite{divfans}
and that way obtain a combinatorial approach 
to the Cox ring, see Theorem~\ref{sec:thm-cox-divfan}.
Similarly to the toric case~\cite{Co}, 
the advantage of the combinatorial treatment is 
that the divisor class group is easily accessible
via the defining data and thus one has a 
simple approach to the grading of the Cox ring.

In Section~\ref{applexam}, we give 
some applications.
The description of the Cox ring 
given in Theorem~\ref{complexity1} allows us 
to apply the language of bunched rings 
presented in~\cite{BeHa2} and~\cite{Ha2} in order to 
investigate complete normal rational 
varieties $X$ with a complexity one 
torus action.
For example, in Corollary~\ref{cor:embedd},
we realize $X$ as an invariant complete 
intersection in a toric variety $X'$,
provided that any two points of $X$ admit a 
common affine neighborhood.
Moreover, in  Corollary~\ref{cor:canondiv},
we obtain explicit descriptions of the cone of 
movable divisor classes and the canonical divisor 
in terms of the divisors $E_k$ and $D_{ij}$.

The first non-trivial examples of complexity 
one torus actions are complete normal rational
$\KK^*$-surfaces $X$.
An important data is the Orlik-Wagreich
graph associated to $X$, which describes 
the intersection theory of a canonical resolution 
$\t{X}$ of $X$, see~\cite{OrWa}.
In Theorem~\ref{OrWa2Cox}, we show how to extract 
the Cox ring of $\t{X}$ from the Orlik-Wagreich
graph, which in turn allows to compute the 
Cox ring of $X$. 
In Theorem~\ref{thm:crgordelps}, we explicitly
compute the  Cox rings 
of all Gorenstein del Pezzo
surfaces of Picard number at most two.

Finally, we consider in Section~\ref{applexam}
projectivizations of 
equivariant vector bundles over 
complete toric varieties.
We explicitly compute the Cox ring for 
the case of rank two bundles and 
for the case of the cotangent bundle 
over a smooth toric variety, 
see Theorems~\ref{sec:cor-bundle}
and~\ref{sec:cor-cotan}.


We would like to thank the referee for carefully reading 
the manuscript and for many helpful remarks and corrections.

\goodbreak

\section{Cox rings and universal torsors}
\label{coxring}

Here we provide basic ingrediences
for the proofs of Theorems~\ref{fingenchar} to~\ref{complexity1},
which also might be of independent interest.
For example, in Proposition~\ref{prop:sepcox}
we determine the Cox ring of a prevariety $X$ in 
terms of that of a separation $X \to Y$
and Proposition~\ref{liftaction} is 
a lifting statement for torus actions 
to the universal
torsor in the case of torsion in the class group.

We work over an algebraically closed field~$\KK$ of
characteristic zero. 
We will not only deal with varieties over $\KK$ 
but more generally with prevarieties, i.e.,
possibly non-separated spaces.
Recall that a ($\KK$-)prevariety is a space $X$ 
with a sheaf $\mathcal{O}_X$ of $\KK$-valued functions
such that $X = X_1 \cup \ldots \cup X_r$ holds
with open subspaces $X_i$, each of which is an affine 
($\KK$-)variety.

In the sequel, $X$ denotes an irreducible normal 
prevariety.
As in the separated case, the group of Weil 
divisors is the free abelian group $\WDiv(X)$ 
generated by all prime divisors, i.e., 
irreducible subvarieties of codimension one.
The divisor class group $\Cl(X)$ is the 
factor group of $\WDiv(X)$ modulo the 
subgroup of principal divisors.
We define the Cox ring of $X$ 
following~\cite[Sec.~2]{Ha2}.
Suppose that 
$\Gamma(X,\mathcal{O}^*) = {\KK^*}$
holds and that the
divisor class group $\Cl(X)$ 
is finitely generated.
Let $\mathfrak{D} \subseteq \WDiv(X)$ 
be a finitely generated subgroup mapping 
onto $\Cl(X)$ and consider the sheaf of 
$\mathfrak{D}$-graded algebras
$$
\mathcal{S}
\ :=  \
\bigoplus_{D \in \mathfrak{D}} \mathcal{S}_D,
\qquad\qquad
\mathcal{S}_D
\ :=  \
\mathcal{O}_X(D),
$$
where multiplication is defined via  
multiplying homogeneous sections as 
rational functions on $X$. 
Let $\mathfrak{D}^0 \subseteq \mathfrak{D}$
be the kernel of $\mathfrak{D} \to \Cl(X)$.
Fix a  {\em shifting family\/}, i.e.,
a family of $\mathcal{O}_X$-module 
isomorphisms 
$\varrho_{D^0} \colon \mathcal{S} \to \mathcal{S}$,
where $D^0 \in \mathfrak{D}^0$, such that 
\begin{itemize}
\item
$\varrho_{D^0}(\mathcal{S}_D) = \mathcal{S}_{D+D^0}$
for all $D \in \mathfrak{D}$, $D^0 \in \mathfrak{D}^0$,
\item
$\varrho_{D^0_1+D^0_2} = \varrho_{D^0_2} \circ \varrho_{D^0_1}$
for all $D^0_1,D^0_2 \in \mathfrak{D}^0$,
\item
$\varrho_{D^0}(fg) = f \varrho_{D^0}(g)$ for all $D^0   \in \mathfrak{D}^0$
and any two homogeneous $f,g$. 
\end{itemize}
Consider the quasicoherent sheaf $\mathcal{I}$
of ideals of $\mathcal{S}$ generated by all 
sections of the form $f - \varrho_{D^0}(f)$, 
where $f$ is homogeneous  
and $D^0$ runs through $\mathfrak{D}^0$.
Then $\mathcal{I}$ is homogeneous with respect 
to the coarsened grading 
$$
\mathcal{S}
\ = \ 
\bigoplus_{[D] \in \Cl(X)} \mathcal{S}_{[D]},
\qquad \qquad 
\mathcal{S}_{[D]}
\ = \ 
\bigoplus_{D' \in D + \mathfrak{D}^0} \mathcal{O}_X(D').
$$
Moreover, it turns out that $\mathcal{I}$ is
a sheaf of radical ideals.
Dividing the $\Cl(X)$-graded $\mathcal{S}$ by the
homogeneous ideal $\mathcal{I}$, 
we obtain a quasicoherent 
sheaf of $\Cl(X)$-graded $\mathcal{O}_X$-algebras, 
the {\em Cox sheaf}: set $\mathcal{R} :=  \mathcal{S}/\mathcal{I}$,
let $\pi \colon \mathcal{S} \to \mathcal{R}$ 
be the projection and define the grading by 
$$
\mathcal{R}
\ = \
\bigoplus_{[D] \in \Cl(X)} \mathcal{R}_{[D]},
\qquad \qquad 
\mathcal{R}_{[D]}
\ = \ 
\pi \left(\mathcal{S}_{[D]} \right).
$$
One can show that, up to isomorphism,
the graded $\mathcal{O}_X$-algebra
$\mathcal{R}$ does not depend on the 
choices of $\mathfrak{D}$ and the shifting family.
We define the {\em Cox ring\/} $\mathcal{R}(X)$ 
of $X$, also called the {\em total coordinate ring\/} 
of $X$, to be the $\Cl(X)$-graded 
algebra of global sections of the Cox sheaf:
$$
\mathcal{R}(X)
\ := \
\Gamma(X,\mathcal{R})
\ \cong \ 
\Gamma(X,\mathcal{S}) /\Gamma(X,\mathcal{I}).
$$

We are ready to perform first computations of
Cox rings. Our aim is to relate the Cox ring 
of a prevariety $X$ to that of a (separated) 
variety arising in a canonical way from $X$.
We say that an open subset $U \subseteq X$ 
is {\em big\/} if the complement $X \setminus U$
is of codimension at least two in $X$.

\begin{definition}
\label{def:goodsep}
By a {\em separation\/} of a prevariety
$X$ we mean a rational map $\varphi \colon X \dashrightarrow Y$
to a (separated) variety $Y$,
which is defined on a big open subset
$U \subseteq X$ and maps $U$ locally isomorphic onto 
a big open subset $V \subseteq Y$.
\end{definition}

\goodbreak

\begin{remark}
\label{rem:goodsep}
Let $\varphi \colon X \dashrightarrow Y$ be a separation.
Then there are big open subsets 
$U \subseteq X$ and $V \subseteq Y$ such that
$\varphi \colon U \to V$ is a
local isomorphism and moreover
there are prime divisors $C_0, \ldots, C_r$ on $V$
such that 
\begin{enumerate}
\item 
$\varphi$ maps 
$U \setminus \varphi^{-1}(C_0 \cup \ldots \cup C_r)$
isomorphically onto  
$V \setminus (C_0 \cup \ldots \cup C_r)$,
\item
Each $\varphi^{-1}(C_i)$ is a 
disjoint union of prime divisors 
$C_{ij}$ of $U$.
\end{enumerate}
\end{remark}

As we will see in Proposition~\ref{prop:sepex},
every prevariety $X$ with finitely generated divisor
class group admits a separation $X \to Y$.
Here comes how the Cox rings $\mathcal{R}(X)$ and 
$\mathcal{R}(Y)$ are related to each other; 
for the sake of a simple notation, we identify 
prime divisors of the big open subsets 
$U \subseteq X$ and $V \subseteq Y$ 
with their closures in $X$ and $Y$ respectively.

\begin{proposition}
\label{prop:sepcox}
Let $\varphi \colon X \dashrightarrow Y$ be a 
separation, $C_0, \ldots, C_r$
prime divisors on $Y$ as in~\ref{rem:goodsep},
and 
$\varphi^{-1}(C_i) = C_{i1} \cup \ldots \cup C_{in_i}$ 
with pairwise disjoint prime divisors 
$C_{ij}$ on $X$. 
Then $\varphi^* \colon \Cl(Y) \to \Cl(X)$
is injective, and we have 
\begin{eqnarray*}
\Cl(X) 
& = & 
\varphi^* \Cl(Y) 
\ \oplus \ 
\bigoplus_{\genfrac{}{}{0pt}{}{0 \le i \le r,}{1 \le j \le n_i-1}}
\ZZ [C_{ij}].
\end{eqnarray*}
If $\Gamma(X,\mathcal{O}^*) = \KK^*$ holds and 
$\Cl(X)$ is finitely generated, then 
there is a canonical injective pullback homomorphism
$\varphi^* \colon \mathcal{R}(Y) \to \mathcal{R}(X)$ 
of Cox rings.
Moreover, with $\deg(T_{ij}) := [C_{ij}]$ 
and $T_i := T_{i1} \cdots T_{in_i}$,
the assignment $T_{ij} \mapsto 1_{C_{ij}}$ defines a
$\Cl(X)$-graded isomorphism
\begin{eqnarray*}
\mathcal{R}(Y)[T_{ij}; \; 0 \le i \le r, \, 1 \le j \le n_i]
\ / \ 
\bangle{T_i - 1_{C_i}; \; 0 \le i \le r}
& \to & 
\mathcal{R}(X).
\end{eqnarray*}
\end{proposition}

\begin{proof}
Since divisor class group and Cox ring do not 
change when passing to big open subsets, we may 
assume $U= X$ and $V = Y$ in the setting of Remark~\ref{rem:goodsep}.
The assertion on the divisor class group follows immediately 
from the facts that the principal divisors of $X$ are precisely 
the pull backs of principal divisors on $Y$ and that the divisor
class group of $X$ is generated by all pullback divisors 
and the classes $[C_{ij}]$, where $0 \le i \le r$ and  
$1 \le j \le n_i-1$.

We turn to the Cox rings. Let $\mathfrak{D}_Y \subseteq \WDiv(Y)$ be a 
finitely generated subgroup containing $C_0, \ldots, C_r$ 
and mapping onto $\Cl(Y)$.
Moreover, let $\mathfrak{D}_X \subseteq \WDiv(X)$ be the subgroup 
generated by $\varphi^*(\mathfrak{D}_Y)$ and the divisors $C_{ij}$, 
where $0 \le i \le r$ and $1 \le j \le n_i-1$; 
note that $C_{in_i} \in \mathfrak{D}_X$ holds. 
Consider the associated graded sheaves
$$ 
\mathcal{S}_Y \ := \ \bigoplus_{E \in \mathfrak{D}_Y} \mathcal{O}_Y(E),
\qquad \qquad
\mathcal{S}_X \ := \ \bigoplus_{D \in \mathfrak{D}_X} \mathcal{O}_X(D).
$$
Then we have a graded injective pullback homomorphism 
$\varphi^* \colon \mathcal{S}_Y \to \mathcal{S}_X$, 
which in turn extends to a homomorphism 
$$ 
\psi \colon
\mathcal{S}_Y[T_{ij}; \; 0 \le i \le r, \, 1 \le j \le n_i]
\ \to \ 
\mathcal{S}_X,
\qquad
T_{ij} 
\ \mapsto \ 
1_{C_{ij}}.
$$

We show that $\psi$ is surjective.
Given a section $h$ of $\mathcal{S}_X$
of degree $D \in \mathfrak{D}_X$,
consider its divisor
$D(h) = \div(h) + D$.
If there occurs a $C_{ij} \in \mathfrak{D}_X$ in $D(h)$, 
then we may divide $h$ in $\mathcal{S}_X$ 
by the corresponding $1_{C_{ij}}$.
Doing this as often as possible, we arrive 
at some section $h'$ of $\mathcal{S}_X$,
homogeneous of some degree $D' \in \mathfrak{D}_X$,
such that $D(h') = \div(h')+D'$
has no components $C_{ij}$.
But then $D'$ is a pullback divisor and 
$h'$ is a pullback section.
This in turn means that $h'$ is a polynomial
over $\varphi^* \mathcal{S}_Y$ and 
the $1_{C_{ij}}$. 

Next, we determine the kernel of $\psi$,
which amounts to determining  the 
relations among the sections 
$s_{ij} := 1_{C_{ij}}$.
Consider two coprime monomials $F, F'$ 
in the $s_{ij}$ and two 
homogeneous pullback sections 
$h, h'$ of $\varphi^*(\mathcal{S}_Y)$.
If $\deg(hF) = \deg(h'F')$ holds in $\mathfrak{D}_X$, 
then the difference $\deg(F') - \deg(F)$
must be a linear combination of some 
$\varphi^*(C_i) \in \mathfrak{D}_X$ and hence   
$F$ and $F'$ are products of some $\varphi^* 1_{C_i}$.
As a consequence, we obtain that any
homogeneous (and hence any)
relation among the $s_{ij}$ 
is generated by 
the relations $T_i - 1_{C_i}$.

Finally, fix a shifting family $\varrho_Y$ for $\mathcal{S}_Y$.
Since $\mathfrak{D}_X^0 = \varphi^*(\mathfrak{D}_Y^0)$ 
holds, the pullback family $\varphi^* \varrho_Y$ extends 
uniquely to a shifting family $\varrho_X$ for $\mathcal{S}_X$.
We have $\mathcal{I}_X = \varphi^*(\mathcal{I}_Y)$ and  
hence obtain a well defined graded pullback homomorphism
$\varphi^* \colon \mathcal{R}(Y) \to \mathcal{R}(X)$,
which is injective, because 
$\varphi^* \colon  \mathfrak{D}_Y /  \mathfrak{D}_Y^0 
\to \mathfrak{D}_X / \mathfrak{D}_X^0$ is so and 
$\varphi^* \colon \mathcal{S}_Y \to \mathcal{S}_X$
is an isomorphism when restricted to homogeneous 
components.
Now one directly verifies that
the above epimorphism $\psi$ induces the desired 
isomorphism.
\end{proof}

We apply this result to compute the 
Cox ring of the prevariety occurring 
as non-separated orbit space for 
torus actions of complexity one.
Consider the projective line  
$\PP_1$, a tuple
$A = (a_0, \ldots, a_r)$
of pairwise different points
$a_i$ on $\PP_1$,
and a tuple  
$\mathfrak{n} = (n_0, \ldots, n_r) \in \ZZ_{\ge 1}^r$,
where $r \ge 1$.
Set
$$ 
X_{ij}
\ := \ 
\PP_1 \setminus \bigcup_{k \ne i} a_k,
\qquad 
0 \le i \le r, 
\qquad
1 \le j \le n_i.
$$
Then, gluing the $X_{ij}$ along the common open 
subset $\PP_1 \setminus \{a_0, \ldots, a_r\}$,
one obtains an irreducible smooth prevariety 
$\PP_1(A,\mathfrak{n})$ of dimension one.
The inclusion maps $X_{ij} \to \PP_1$ glue together 
to a morphism $\pi \colon \PP_1(A,\mathfrak{n}) \to \PP_1$,
which is a separation.
Writing $a_{ij}$ for the point in $\PP_1(A,\mathfrak{n})$
stemming from $a_i \in X_{ij}$, we obtain the fiber
over a point $a \in \PP_1$ as 
\begin{eqnarray*}
\pi^{-1}(a)
& = &
\begin{cases}
\{a_{i1}, \ldots, a_{in_i}\} & a = a_i \text{ for some } 0 \le i \le r,
\\
\{a\}                        & a \ne a_i \text{ for all } 0 \le i \le r.
\end{cases}
\end{eqnarray*}
For every $0 \le i \le r$,
define a monomial
$T_i := T_{i1} \cdots T_{in_i}$
in the polynomial ring 
$\KK[T_{ij}; \; 0 \le i \le r, \; 1 \le j \le n_i]$.
Moreover, for every $a_i \in \PP_1$ 
fix a presentation 
$a_i = [b_i,c_i]$ with $b_i,c_i \in \KK$
and for every $0 \le i \le r-2$ set 
$k = j+1 = i+2$ 
and define a trinomial
\begin{eqnarray*}
g_i 
& := &
(c_kb_j - c_jb_k)T_i
\ + \ 
(c_ib_k - c_kb_i)T_j
\ + \ 
(c_jb_i - c_ib_j)T_k.
\end{eqnarray*}

\begin{proposition}
\label{P1AnCox}
The divisor class group of 
$\PP_1(A,\mathfrak{n})$ is free of rank 
$n_0 +\ldots + n_r -r$ and
there is a decomposition
\begin{eqnarray*}
\Cl(\PP_1(A,\mathfrak{n}))
& = &
\bigoplus_{j=1}^{n_0} \ZZ \mal [a_{0j}] 
\ \oplus \
\bigoplus_{i=1}^r
\left( \bigoplus_{j=1}^{n_i-1} \ZZ \mal [a_{ij}] 
\right).
\end{eqnarray*}
Moreover, in terms of the above data
and with $\deg(T_{ij}) := [a_{ij}]$,
the Cox ring of 
$\PP_1(A,\mathfrak{n})$ is given as
\begin{eqnarray*}
\mathcal{R}(\PP_1(A,\mathfrak{n}))
&  \cong &
\KK[T_{ij}; \; 0 \le i \le r, \; 1 \le j \le n_i]
\  / \ 
\bangle{g_i; \; 0 \le i \le r-2}.
\end{eqnarray*}
\end{proposition}

\begin{proof}
The statement on the divisor class group
is clear.
The description of the 
Cox ring follows from Proposition~\ref{prop:sepcox}
and the fact that the Cox ring $\mathcal{R}(\PP_1)$ 
of the projective line is generated by the 
canonical sections $s_i := 1_{a_i}$ and
has the relations 
\begin{eqnarray*}
(c_kb_j - c_jb_k)s_i
\ + \ 
(c_ib_k - c_kb_i)s_j
\ + \ 
(c_jb_i - c_ib_j)s_k
& = &
0,
\end{eqnarray*}
where
$0 \le i \le r-2$,
and $k = j + 1 = i+2$;
note that the dependence of these relations
on the choice of the $b_i,c_i$ reflects
the choice of a shifting family.
\end{proof}

Now we discuss some geometric 
aspects of the Cox ring.
As before, let $X$ be a normal prevariety with 
$\Gamma(X,\mathcal{O}^*) = \KK^*$
and finitely generated divisor class 
group, 
and let $\mathcal{R}$ be a Cox sheaf.
Suppose that $\mathcal{R}$ is locally of finite type;
this holds for example if $X$ is locally factorial
or if $\mathcal{R}(X)$ is finitely generated.
Then we may consider the relative spectrum 
\begin{eqnarray*}
\rq{X}
& := & 
\Spec_X(\mathcal{R}).
\end{eqnarray*}
The $\Cl(X)$-grading of the sheaf $\mathcal{R}$ 
of $\mathcal{O}_X$-algebras
defines an action of the diagonalizable group
$H_X := \Spec \, \KK[\Cl(X)]$ on $\rq{X}$,  
and the canonical morphism $p \colon \rq{X} \to X$
is a good quotient, i.e., it is an $H_X$-invariant
affine morphism satisfying
\begin{eqnarray*}
\mathcal{O}_X 
& = & (p_*\mathcal{O}_{\rq X})^{H_X}.
\end{eqnarray*}
We call $p \colon \rq{X} \to X$ the {\em universal 
torsor\/} associated to $\mathcal{R}$.
If the Cox ring $\mathcal{R}(X)$ 
is finitely generated, then we define the 
{\em total coordinate space\/}
of $X$ to be the affine variety 
$\b{X} = \Spec(\mathcal{R}(X))$
together with the $H_X$-action defined 
by the $\Cl(X)$-grading of $\mathcal{R}(X)$.

As usual, we say that a Weil divisor 
$\sum a_D D$, where $D$ runs through 
the irreducible 
hypersurfaces, on a prevariety $Y$ with an 
action of a group $G$
is {\em $G$-invariant\/} if $a_D = a_{g \cdot D}$ 
holds for all $g \in G$. 
We say that $Y$ is {\em $G$-factorial\/} if 
every $G$-invariant divisor on $G$ 
is principal.
Moreover,  we say that a prevariety $Y$ is of
{\em affine intersection\/} if for any two affine 
open subsets $V,V' \subseteq Y$ the intersection 
$V \cap V'$ is again affine.

\begin{proposition}
\label{gencoxprops}
Let $X$ be an irreducible smooth prevariety 
of affine intersection with 
$\Gamma(X, \mathcal{O}^*) = \KK^*$ and
finitely generated divisor class group.
Let $\mathcal{R}$ be a Cox sheaf and
denote by $p \colon \rq{X} \to X$ 
the associated universal torsor.
\begin{enumerate}
\item
$\rq{X}$ is a normal quasiaffine 
variety, and every homogeneous invertible 
function on $\rq{X}$ is constant.
If $\Gamma(X, \mathcal{O}) = \KK$ holds
or $\Cl(X)$ is free, 
then even every invertible function 
on $\rq{X}$ is constant.
\item
The action of $H_X$ on $\rq{X}$ is free
and $\rq{X}$ is $H_X$-factorial.
If $\Cl(X)$ is free, 
then $\rq{X}$ is even factorial.
\end{enumerate}
\end{proposition}

\begin{proof}
Normality of $\rq{X}$ follows 
from~\cite[Lemma~3.10]{BeHa1}.
Since $X$ is of affine intersection, 
it can be covered by open affine subsets,
the complements of which are of pure codimension 
one.
Together with smoothness this implies that 
$X$ is divisorial in the sense 
of~\cite[Sec.~4]{BeHa1}.
Thus, we infer from~\cite[Prop.~6.3]{BeHa1}
that $\rq{X}$ a quasiaffine variety.
The fact that every homogeneous invertible 
function on $\rq{X}$ is constant is seen 
as in ~\cite[Prop.~2.2~(i)]{Ha2}.
Moreover, \cite[Thm.~7.3]{BeHa1} tells
us that every invertible 
function on $\rq{X}$ is constant if
we have $\Gamma(X,\mathcal{O}) = \KK$.
For Assertion~(ii), we can proceed 
exactly as in the proof 
of~\cite[Prop.~2.2~(iv)]{Ha2}.
\end{proof}

\begin{proposition}
\label{liftaction}
Let $X$ be an irreducible smooth prevariety 
of affine intersection with 
$\Gamma(X, \mathcal{O}^*) = \KK^*$ and
finitely generated divisor class group.
Let $\mathcal{R}$ be a Cox sheaf on~$X$
and $p \colon \rq{X} \to X$ 
the associated universal torsor.
Assume that $T \times X \to X$ is
an effective algebraic torus action.

\begin{enumerate}
\item
There is a $T$-action on $\rq{X}$ 
and an
epimorphism $\varepsilon \colon T \to T$ 
such that for all $h \in H_X$, 
$t \in T$ and $z \in \rq{X}$ one has
$$
t \mal h \mal z \ = \ h \mal t \mal z,
\qquad
p(t \mal z) \ = \ \varepsilon(t) \mal p(z).
$$
If the divisor class group $\Cl(X)$ is free, 
then one may take the homomorphism 
 $\varepsilon \colon T \to T$ to be the 
identity.
\item
Let $T \times H_X$ act on $\rq{X}$ as 
in~(i),
let $G' \subseteq T \times H_X$ be the trivially
acting subgroup and consider the induced
effective action of $G := (T \times H_X)/G'$ 
on $\rq{X}$. 
Then for any $z \in \rq{X}$, there is an 
isomorphism of isotropy groups $G_z \cong T_{p(z)}$.
\end{enumerate}
\end{proposition}

\begin{proof}
We prove~(i).
Take a group $\mathfrak{D} \subseteq \WDiv(X)$ of 
Weil divisors mapping onto the divisor class group, 
and let $D_1, \ldots, D_r \in \WDiv(X)$ be a basis 
of $\mathfrak{D}$ such that the kernel 
$\mathfrak{D}_0 \subseteq \mathfrak{D}$ of 
$\mathfrak{D} \to \Cl(X)$ has a basis
of the form  $a_iD_i$, 
where $1 \le i \le s$ with some $s \le r$.

For every $D_i$ choose a $T$-linearization,
and via tensoring these linearizations, define
a $T$-linearization of the whole group $\mathfrak{D}$,
compare~\cite[Sec.~1]{Ha1}. 
Note that the $T$-linearization of the trivial
divisor $a_iD_i$ is given by a character
$\chi_i$.
Set $b := a_1 \cdots a_s$ and consider the 
epimorphism $\varepsilon \colon T \to T$,
$t \mapsto t^b$.
Then we have a new $T$-action  
$$
T \times X \ \to \ X,
\qquad
(t,x) \ \mapsto \ \varepsilon(t) \mal x.
$$
The divisors $D \in \mathfrak{D}$
are as well linearized with respect to this 
new action.
Twisting each $T$-linearization 
of $D_i$ with $\chi_i^{-b/a_i}$,
we achieve that each $a_iD_i$ is trivially
$T$-linearized with respect to the 
new $T$-action on $X$.
Thus, we may choose $T$-equivariant
isomorphisms 
$\varrho_{i} \colon \mathcal{O}_X \to \mathcal{O}_X(a_iD_i)$.

Let $\mathcal{S}$ denote the $\mathfrak{D}$-graded 
sheaf defined by $\mathfrak{D}$. 
Using the isomorphisms $\varrho_i$, we construct a 
$T$-equivariant shifting family:
for $D^0 = b_1a_1D_1 + \cdots + b_sa_s D_s$, 
define a $T$-equivariant isomorphism
$\varrho_{D^0} \colon \mathcal{S} \to \mathcal{S}$ 
by sending a $\mathfrak{D}$-homogeneous $f$ to
\begin{eqnarray*}
\varrho_{D^0}(f) 
& := & 
\varrho_{1}(1)^{b_1}
\cdots
\varrho_{s}(1)^{b_s} 
f.
\end{eqnarray*}
The ideal $\mathcal{I}$ of $\mathcal{S}$ associated
to this shifting family is $T$-homogeneous.
This means that the $T$-action on $\Spec_X \mathcal{S}$ 
defined by the $T$-linearization of $\mathfrak{D}$
leaves $\rq{X}$ invariant. 
By construction, the torsor $p \colon \rq{X} \to X$
is $T$-equivariant, when we take the 
new $T$-action on $X$.

We turn to~(ii). 
Let $\varepsilon \colon T \to T$ be as in~(i).
A first step is to show 
that for any given point $z \in \rq{X}$, 
the kernel of ineffectivity 
$G' \subseteq T \times H_X$ 
can be written as
\begin{eqnarray*}
G'
&  =  &
\{(t,h) \in (T \times H_X)_z; \; \varepsilon(t) = 1 \}.
\end{eqnarray*}

In order to verify the inclusion ``$\subseteq$'',
let $(t,h) \in G'$ be given.
Then $(t,h) \mal z' = z'$ holds for every 
point $z' \in \rq{X}$.
In particular, $(t,h)$ belongs to $(T \times H_X)_z$.
Moreover, we obtain $\varepsilon(t) \mal p(z') = p(z')$
for every $z' \in \rq{X}$.
Since $p \colon \rq{X} \to X$ is surjective and
$T$ acts effectively on $X$, this implies
$\varepsilon(t) = 1$.

For checking the inclusion ``$\supseteq$'',
consider $(t,h) \in (T \times H_X)_z$ with 
$\varepsilon(t) = 1$.
Then, for every $z' \in \rq{X}$, we have 
$p((t,h) \mal z') = p(z')$.
Consequently $t \mal z' = h(t,z') \mal z'$ 
holds with a uniquely determined $h(t,z') \in H_X$. 
Consider the assignment
$$ 
\eta \colon \rq{X} \ \to \ H_X, 
\qquad 
z' \ \mapsto \ h(t,z').
$$
Since $H_X$ acts freely we may choose for 
any $z'$ homogeneous functions $f_1, \ldots, f_r$,
defined near $z'$ with $f_i(z') = 1$ such
that their weights $\chi_1, \ldots, \chi_r$ 
form a basis of the character group of $H_X$.
Then, near $z'$, we have a commutative diagram
$$ 
\xymatrix{
&& 
{\rq{X}}
\ar[dll]_{\eta}
\ar[drr]^{\qquad z' \ \mapsto \ (f_1(t \cdot z'), \ldots, f_r(t \cdot z'))}
&&
\\
H_X
\ar[rrrr]^{\cong}_{h' \ \mapsto \ (\chi_1(h'), \ldots, \chi_r(h'))}
&&&&
{(\KK^*)^r}
}
$$
Consequently, the map $\eta$ is 
is a morphism.
Moreover, pulling back characters of $H_X$ via
$\eta$ gives invertible $H_X$-homogeneous 
functions on $\rq{X}$, which  
by Proposition~\ref{gencoxprops}~(i) 
are constant. 
Thus, $\eta$ is constant, which means 
that $h(t) := h(t,z')$ does not depend
on $z'$. 
By construction, $(t,h(t)^{-1})$
belongs to $G'$. 
Moreover, $t \mal z = h^{-1} \mal z$ 
and freeness of the $H_X$-action give
$h(t) = h^{-1}$.
This implies $(t,h) \in G'$.

We are ready to prove the assertion.
Note that $(t,h) \mapsto \varepsilon(t)$ 
defines a homomorphism 
$\beta \colon (T \times H_X)_z \to T_{p(z)}$.
We claim that $\beta$ is surjective.
Given $t \in T_{p(z)}$, 
choose $t' \in T$ with $\varepsilon(t') = t$.
Then we have
$$ 
p(t' \mal z)
\ = \ 
\varepsilon(t') \mal p(z) 
\ = \ 
t \mal p(z) 
\ = \
p(z).
$$
Consequently, $t' \mal z = h \mal z$ holds
for some $h \in H_X$. Thus, 
$(t',h^{-1}) \in (T \times H_X)_z$
is mapped by $\beta$ to $t \in T_{p(z)}$.   
By the first step, the kernel of $\beta$ 
is just $G'$. This gives a commutative
diagram
$$ 
\xymatrix{
(T \times H_X)_z
\ar[rr]^{\beta}
\ar[dr]_{/G'}
&&
T_{p(z)}
\\
&
G_z
\ar[ur]_{\cong}
&
}
$$
\end{proof}

In the sequel, we mean by a 
{\em universal torsor\/} for $X$
more generally any good quotient 
$q \colon \mathcal{X} \to X$ for an action 
of $H_X$ on a variety $ \mathcal{X}$ 
such that there is an equivariant 
isomorphism 
$\imath \colon  \mathcal{X} \to \rq{X}$ with 
$q = p \circ \imath$.

\begin{proposition}
\label{bigfreequot}
Let $\mathcal{X}$ be a normal quasiaffine 
variety with a free action of 
a diagonalizable group $H$.
If every invertible function on
$\mathcal{X}$ is constant 
and $\mathcal{X}$ is $H$-factorial, 
then the quotient $q \colon \mathcal{X} \to X$
is a universal torsor for 
$X := \mathcal{X} / H$.
\end{proposition}

\begin{proof}
We have $H = \Spec \, \KK[K]$ with the character
lattice $K$ of $H$.
A first step is to provide an isomorphism 
$K \to \Cl(X)$.
Cover $\mathcal{X}$
by $H$-invariant affine open subsets $W_{j}$ 
such that, for every $w \in K$ and every $j$,
there is a $w$-homogeneous 
$h_{w,j} \in \Gamma(W_{j},\mathcal{O}^{*})$.
Moreover, for every $w \in K$, fix a $w$-homogeneous 
$h_{w} \in \KK(\mathcal{X})^{*}$.
Then the $H$-invariant local equations $h_w/h_{w,j}$ 
define a Weil divisor $D(h_w)$ on $X$ 
satisfying
$$
D(h_w) \ = \ \div(h_{w}/h_{w,j})
\text{ on } q(W_j),
\qquad \qquad 
q^*(D(h_w)) \ = \ \div(h_w).
$$
We claim that the assignment $w \mapsto D(h_w)$ 
induces an isomorphism 
from $K$ onto $\Cl(X)$, not depending on the choice of 
$h_{w}$:
$$
K \ \to \ \Cl(X),
\qquad
w \ \mapsto \ \b{D}(w) := [D(h_w)].
$$
To see that the class $\b{D}(w)$ does not 
depend on the choice of $h_w$, consider a further 
$w$-homogeneous $g_w \in \KK(X)^*$. 
Then $f := h_w/g_w$ is an
invariant rational function descending 
to $X$, where we obtain
\begin{eqnarray*}
D(h_w) \ - \ D(g_w) &  = & \div(f).
\end{eqnarray*}
Thus, $K \to \Cl(X)$ is a well defined 
homomorphism.
To verify injectivity, let $D(h_w) = \div(f)$ 
for some $f \in \KK(X)^*$. 
Then we obtain $\div(h_w) = \div(q^*(f))$.
Thus, $h_w/q^*(f)$ is an invertible homogeneous
function on $\mathcal{X}$ 
and hence is constant.
This implies $w = 0$.
For surjectivity, let any 
$D \in \WDiv(X)$ be given.
Then $q^*(D)$ is an $H$-invariant divisor
on $\mathcal{X}$ and hence we have 
$q^*(D) = \div(h)$ with some rational
function $h$ on $\mathcal{X}$, which
is homogeneous of some degree
$w$. This means $D = D(h)$.

Now, choose a group 
$\mathfrak{D} \subseteq \WDiv(X)$ of 
Weil divisors mapping onto the divisor
 class group $\Cl(X)$, 
and let $D_1, \ldots, D_r \in \WDiv(X)$ be a basis 
of $\mathfrak{D}$ such that the kernel 
$\mathfrak{D}_0 \subseteq \mathfrak{D}$ of 
$\mathfrak{D} \to \Cl(X)$ is generated by 
multiples $a_iD_i$, 
where $1 \le i \le s$ with some $s \le r$.
Set
$$
\mathcal{S}
\ := \
\bigoplus_{D \in \mathfrak{D}}
\mathcal{S}_D,
\qquad \qquad
\mathcal{S}_D
\ := \ 
\mathcal{O}_X(D).
$$
By the preceding consideration,
we may assume $D_i = D(h_{w_i})$ for $1 \le i \le r$.
Then, for $D = b_1D_1 + \ldots + b_rD_r$,
we have $D = D(h_w)$ with 
$h_w = h_{w_1}^{b_1} \cdots  h_{w_r}^{b_r}$,
where $w = [D]$ is the class of 
$D = D(h_w)$ in $K = \Cl(X)$.
For any open $U \subseteq X$, 
we have an
isomorphism of $\KK$-vector spaces
$$ 
\Phi_{U,D}
\colon 
\Gamma(U,\mathcal{O}(D(h_w))) 
\ \to \  
\Gamma(q^{-1}(U), \mathcal{O})_w,
\qquad
g \ \mapsto \ q^*(g) h_w.
$$
In fact, on each $U_j := q(W_j) \cap U$ the 
section $g$ is  given as $g = g_j'h_j/h_w$ 
with a regular function 
$g_j' \in \mathcal{O}(U_j)$.
Consequently, the function $q^*(g) h_w$ is 
regular on $q^{-1}(U)$.
In particular, the assignment is a well defined
homomorphism.
Moreover, $f \mapsto f/h_w$ defines an inverse 
homomorphism.

The $\Phi_{U,D}$ fit together to an epimorphism
of sheaves  
$\Phi \colon 
\mathcal{S} \to q_*(\mathcal{O}_{\mathcal{X}})$.
We claim that the kernel 
$\mathcal{I}$ of $\Phi$ is the ideal
of a shifting family $\varrho$.
Indeed, for any $D^0 \in \mathfrak{D}^0$, consider 
$$ 
h^0 
\ : = \ 
\Phi_{X,D^0}^{-1}(\Phi_{X,0}(1))
\ \in \
\Gamma(X, \mathcal{S}_{D^0}).
$$
Then 
$\varrho_{D_0} \colon \mathcal{S}_{D} \to \mathcal{S}_{D+D^0}$,
$g \mapsto h^0g$ is as wanted.
Thus, we obtain an induced isomorphism 
$\mathcal{R} \to q_* \mathcal{O}_{\mathcal{X}}$,
where $\mathcal{R} = \mathcal{S} / \mathcal{I}$ 
is the associated Cox sheaf. 
This in turn defines the desired isomorphism
$ \mathcal{X} \to \rq{X}$.
\end{proof}

\section{Proof of Theorems~\ref{fingenchar},  \ref{fingenchar2}  
and~\ref{complexity1}}

We begin with a couple of elementary
observations.
Let a diagonalizable group $G$ 
act effectively on a normal 
quasiaffine variety~$U$.
Recall that a function 
$f \in \Gamma(U, \mathcal{O})$ 
is said to be $G$-homogeneous 
of weight $\chi \in \Chi(G)$
if one has $f(g \mal x) = \chi(g) f(x)$ 
for all $g \in G$ and $x \in X$.

\begin{lemma}
\label{lem:fixinv}
If there is a $G$-fixed point $x \in U$, 
then every $G$-homogeneous function
$f \in \Gamma(U, \mathcal{O})$ with 
$f(x) \ne 0$ is $G$-invariant.
\end{lemma}

\begin{proof}
Let $\chi \in \Chi(G)$ be the 
weight of $f \in \Gamma(U, \mathcal{O})$.
Then, for every $g \in G$, we have
$f(x) = \chi(g) f(x)$, which implies
$\chi(g) = 1$.
Thus, $\chi$ is the trivial character.
\end{proof}

By a $G$-prime divisor on $U$ we mean 
a $G$-invariant Weil divisor $\sum a_D D$,
where $D$ runs through the prime divisors, 
we always have $a_D \in \{0,1\}$ and $G$ 
permutes transitively the $D$ with $a_D = 1$.
Let $B_1, \ldots, B_m \subseteq U$ 
be $G$-prime divisors and 
suppose that there are homogeneous
functions 
$f_1, \ldots, f_m \in \Gamma(U, \mathcal{O})$ 
that satisfy $\div(f_i) = B_i$.
Let $\chi_i \in \Chi(G)$ be the weight of $f_i$.

\begin{lemma}
\label{lem:Gitorus}
For $i = 1, \ldots, m$, let $G_i \subseteq G$ 
be the generic isotropy group of $B_i$
and set $G_0 := G_1 \cdots G_m \subseteq G$.
\begin{enumerate}
\item
The restriction of $\chi_i$ to $G_i$ generates
the character group $\Chi(G_i)$.
\item
For any two $i,j$ with $j \ne i$, 
the function $f_i$ is $G_j$-invariant.
\item
The group  $G_0$ is 
isomorphic to the direct product of the 
$G_i \subseteq G$.
\item
$\Gamma(U, \mathcal{O})$ 
is generated by $f_1, \ldots, f_m$
and the $G_0$-invariant functions of 
$U$.
\end{enumerate}
\end{lemma}

\begin{proof}
Choose $G$-invariant affine open subsets
$U_i \subseteq U$ such that
$A_i := U_i \cap B_i$ is non-empty
and $U_i \cap B_j$ is empty for every
$j \ne i$.

To prove~(i), let 
$\xi_i \in \Chi(G_i)$ be given.
Then $\xi_i$ is the restriction of 
some $\eta_i \in \Chi(G)$.
Let $V_i \subseteq U_i$ be a
$G$-invariant affine open subset 
on which $G$ acts freely,
and choose a non-trivial 
$G$-homogeneous function $h_i$ 
of weight $\eta_i$ on $V_i$.
Suitably shrinking $U_i$, we achieve 
that $h_i$ is regular without
zeroes on $U_i \setminus A_i$.
Then, on $U_i$, the divisor 
$\div(h_i)$ is a multiple of the 
$G$-prime divisor $A_i = \div(f_i)$ and hence 
$h_i = a_if_i^k$ holds with a $G$-homogeneous 
invertible function $a_i$ on $U_i$
and some $k \in \ZZ$.
By Lemma~\ref{lem:fixinv},
the function $a_i$ is $G_i$-invariant.
We conclude $\eta_i = k \chi_i$ on $G_i$.

Assertion~(ii) is clear by
Lemma~\ref{lem:fixinv}.
To obtain~(iii), it suffices to show
that $\chi_i$ is trivial on $G_j$ 
for any two $i,j$ with $j \ne i$.
But, according to~(ii), we have 
$f_i = \chi_i(g) f_i$
for every $g \in G_j$,
which gives the claim.
Finally, we verify~(iv).
Given a $G$-homogeneous function
$f \in \Gamma(U, \mathcal{O})$,
we may write $f = f'f_1^{\nu_1} \cdots f_m^{\nu_m}$ 
with $\nu_i \in \ZZ_{\ge 0}$ and a regular 
function $f'$ on $U$, 
which is homogeneous with respect to some 
weight $\chi' \in \Chi(G)$
and has order zero along each $G$-prime 
divisor $B_i$.
By Lemma~\ref{lem:fixinv}, the 
function $f'$ is invariant under every $G_i$
and thus under $G_0$.
\end{proof}

Now we specialize to the case that
$B_1, \ldots, B_m \subseteq U$ 
are precisely the $G$-prime divisors of $U$, 
which are contained in $U \setminus U_0$,
where we set
$$ 
U_0
\ := \ 
\{z \in U; \; \dim(G_z) = 0\}
\ \subseteq \ 
U.
$$

\goodbreak

\begin{proposition}
\label{prop:1isosplit}
For $i = 1, \ldots, m$, let $G_i \subseteq G$ 
be the generic isotropy group of $B_i$
and set $G_0 := G_1 \cdots G_m \subseteq G$.
\begin{enumerate}
\item
Each $G_i$ is a one-dimensional torus.
Moreover, there is a non-empty $G_0$-
invariant open subset $U' \subseteq U$
such that each $B_i$ intersects the closure
of any orbit $G_i \mal z \subseteq U'$.
\item
The $G_0$-action on $U_0$ is free,
admits a geometric quotient 
$\lambda_0 \colon U_0 \to V_0$
and the isotropy groups of 
the induced action of
$H := G/G_0$ on $V_0$ satisfy
$H_{\lambda_0(x)} \cong G_x$
for every $x \in U_0$.
\item
$V_0$ is quasiaffine and, moreover,
if $U$ is $G$-factorial
(admits only constant 
invertible functions), then $V_0$ 
is $H$-factorial 
(admits only constant 
invertible functions).
\item
Every $G_0$-invariant rational function of $U$
has neither poles nor zeroes along outside 
$U_0$.
Moreover, there is an isomorphism
$$
\Gamma(U_0, \mathcal{O})^{G_0}[S_1, \ldots, S_m]
\ \to \ 
\Gamma(U, \mathcal{O}),
\qquad
S_i \ \mapsto \ f_i.
$$
\end{enumerate}
\end{proposition}

\begin{proof}
We prove~(i).
By Lemma~\ref{lem:Gitorus}~(i), every $G_i$ 
is a one-dimensional torus.
To proceed,
take any $G_0$-equivariant affine closure
$U \subseteq \b{U}$ and consider the 
quotient
$\lambda_i \colon \b{U} \to \b{U} \quot G_i$.
It maps the fixed point set of the $G_i$-action 
isomorphically onto its image in the quotient
space $\b{U} \quot G_i$.
Since $\b{U} \quot G_i$ is irreducible and 
of dimension at most $\dim(\b{U})-1$, 
we obtain $\lambda_i(\b{B}_i) = \b{U} \quot G_i$
for the closure $\b{B}_i$ of $B_i$ in $\b{U}$.
It follows that $\b{B}_i$ is irreducible, 
equals the whole fixed point set of $G_i$ in $\b{U}$
and any $G_i$-orbit of $\b{U}$ has a point 
of $\b{B}_i$ in its closure.

We turn to~(ii).
Since none of the $f_i$ has a zero 
inside $U_0$, we infer from 
Lemma~\ref{lem:Gitorus} that $G_0$ acts freely 
on $U_0$.
In particular, the action of $G_0$ on $U_0$
admits a geometric quotient
$\lambda_0 \colon U_0 \to V_0$ 
with a prevariety $V_0$.
The statement on the isotropy groups of 
the $H$-action
on $V_0$ is obvious.

We prove the statements made in~(iii) and~(iv).
Denoting by $\TT^m$ the standard $m$-torus
$(\KK^*)^m$, we have a well defined morphism of 
normal prevarieties
$$ 
\varphi \colon
U_0 \ \to \ V_0 \times \TT^m,
\qquad
x \ \mapsto \ (\lambda_0(x),f_1(x), \ldots, f_m(x)).
$$
According to Lemma~\ref{lem:Gitorus}, the 
weight $\chi_i$ of $f_i$ generates the character 
group of $G_i$ for $i = 1, \ldots, m$.
Using this and the fact that $G_0$ is the direct 
product of $G_1, \ldots, G_m$, 
we conclude that $\varphi$ is bijective
and thus an isomorphism. 
In particular, we conclude that $V_0$ is quasiaffine,
because $U$ and hence $U_0$ is so.

Now, endow $V_0$ with the induced 
action of $H = G/G_0$ and $\TT^m$ with the diagonal
$G_0$-action given by the weights 
$\chi_1, \ldots, \chi_m$ of $f_1, \ldots, f_m$.
Then $\varphi$ becomes $G$-equivariant, 
where $G$ acts via the splitting
$G = H \times G_0$ on $V_0 \times \TT^m$.
Using this, we see that $G$-factoriality 
of $U_0$ implies
$H$-factoriality of $V_0$.

We show now that every $G_0$-invariant 
rational function $f \in \KK(U_0)$ 
has neither zeroes nor poles outside $U_0$.
Recall that $U \setminus U_0$ is the 
union of the zero sets $B_i$ of $f_i$, which in
turn are the fixed point sets of the $G_i$-actions
on $U$.
Since the general orbit $G_0 \mal x \subseteq U$
has a point $x_i \in B_i$ in its 
closure, we see that $f$ has neither poles nor 
zeroes along the $B_i$.
In particular, if $f$ is regular on $U_0$ 
then it is so on the whole~$U$.
As a consequence, we obtain that every 
invertible function on $V_0$ 
is constant provided the same holds for $U$.

Finally, according to Lemma~\ref{lem:Gitorus}~(iv),
the homomorphism of~(iv) is surjective. 
Moreover, since the weights $\chi_1, \ldots, \chi_m$ 
of the $f_1, \ldots, f_m$ are a basis of the 
character group of $G_0$,
there are no relations among the $f_i$.
\end{proof}

Let a diagonalizable group $H$ act effectively 
with at most finite isotropy groups on a 
quasiaffine variety $V$.
Suppose that $V$ is $H$-factorial
and admits only constant  
invertible functions.
Denote by $C_1, \ldots, C_n \subseteq V$ 
those $H$-prime divisors of $V$, on
which $H$ acts with a non trivial 
generic isotropy 
group $H_j$ of order $l_j > 1$
and let $g_1, \ldots, g_n$ be homogeneous 
functions on $V$ with $\div(g_j) = C_j$.

\begin{proposition}
\label{lem:quaff2free}
Consider the action of 
$H_0 := H_1 \cdots H_n \subseteq H$ on 
$V$, and let $\kappa \colon V \to W$ 
be the associated quotient.
\begin{enumerate}
\item
$W$ is a quasiaffine variety 
with an effective 
induced action of $H/H_0$,
and
$h \mapsto\kappa^*(h)$
and 
$T_j \mapsto g_j$
define an isomorphism
$$
\Gamma(W,\mathcal{O})[T_1, \ldots, T_n] 
/ 
\bangle{T_j^{l_j} - g_j^{l_j}; \; j = 1, \ldots, n} 
\ \to \ 
\Gamma(V,\mathcal{O}).
$$
\item
$W$ admits an $(H/H_0)$-factorial
big open subset $W_0 \subseteq W$ 
such that $H/H_0$ acts freely on $W_0$
and $W_0$ has only constant invertible
functions.
\end{enumerate}
\end{proposition}

\begin{proof}
We prove~(i).
By Lemma~\ref{lem:Gitorus}~(i),
every $H_j$ is cyclic.
Moreover,
Lemma~\ref{lem:Gitorus}~(iv)
tells us 
that there is an epimorphism
$$ 
\Gamma(W,\mathcal{O})[T_1, \ldots, T_n] 
\ \to \ 
\Gamma(V,\mathcal{O}),
\qquad
h \ \mapsto \ \kappa^*(h),
\quad
T_j \ \mapsto \ g_j.
$$
{From} Lemma~\ref{lem:Gitorus}~(iii) we infer 
that
$H_0 \subseteq H$ is the direct product of 
$H_1, \ldots, H_n \subseteq H$.
Thus, the quotient 
$\kappa \colon V \to W$ can as well be 
obtained by dividing stepwise by effective 
actions of the~$H_j$.
Using this,
one directly checks that the kernel of this
epimorphism is the ideal generated by 
$T_j^{l_j} - g_j^{l_j}$, where $1 \le j \le n$.

We turn to~(ii).
Note that $W$ admits only constant 
invertible functions.
Let $V_0 \subseteq V$ 
denote the subset consisting of 
all points $y \in V$ that have 
either trivial isotropy group $H_{0,y}$
or belong to some $C_j$ and have 
isotropy group $H_{0,y} = H_j$.
Note that $V_0 \subseteq V$ is 
big, $H$-invariant and open.
Set $W_0 := \kappa(V_0)$.
Then $W_0 \subseteq W$ is 
big and the restriction 
$\kappa \colon V_0 \to W_0$  
is a quotient for the action of $H_0$.
By construction, $H/H_0$ acts freely 
on~$W_0$.

We show that $W_0$ is 
$H/H_0$-factorial.
Since $V_0$ and $W_0$
are normal, there is a smooth 
$(H/H_0)$-invariant big open subset 
$W_1 \subseteq W_0$ 
such that $V_1 := \kappa^{-1}(W_1)$
is also smooth and big in 
$V_0$.
We have to show that 
every $(H/H_0)$-linearizable bundle
on $W_1$ is trivial.
According to~\cite[Cor.~5.3]{KKV},
we have two exact sequences fitting 
into the following diagram
$$ 
\xymatrix{
& 
& 
1 
\ar[d]
&
\\
& 
& 
{\Pic}(W_1) 
\ar[d]^{\kappa^*}
&
\\
1 \ar[r]
&
{\Chi(H_0)} \ar[r]^{\alpha}
&
{\Pic_{H_0}}(V_1) \ar[r]^{\beta}
\ar[d]^{\delta}
&
{\Pic}(V_1)
\\
& 
& 
{\prod_{i=1}^n} \Chi(H_i) 
&
}
$$
where the isotropy groups $H_1, \ldots, H_n$ 
generate $H_0$ and hence $\beta \circ \kappa^*$ 
and $\delta \circ \alpha$ are injective.
Given an $(H/H_0)$-linearized 
bundle $L$ on $W_1$, the pullback
$\kappa^*(L)$ is trivial by assumption,
which means $\beta(\kappa^*(L)) = 1$.
Consequently, $L$ is trivial.
\end{proof}

\begin{proof}{Proof of Theorem~\ref{fingenchar}}
We prove the statement more generally for the 
cases that 
$\Gamma(X, \mathcal{O}) = \KK$ holds
or that 
$\Gamma(X, \mathcal{O}^*) = \KK^*$ holds 
and $\Cl(X)$ is free.
Since the Cox ring of $X$ and that of its 
smooth locus coincide, we may assume that 
$X$ is smooth.
Consider the universal torsor 
$p \colon \rq{X} \to X$.
By Proposition~\ref{gencoxprops},
this is a geometric quotient for a free action 
of the diagonalizable group
$H_X := \Spec \, \KK[\Cl(X)]$
on $\rq{X}$
and $\rq{X}$ has only constant globally invertible 
functions.
Fix a lifting of the $T$-action
to $\rq{X}$ as in Proposition~\ref{liftaction}~(i)
and, as in Proposition~\ref{liftaction}~(ii),
let $G$ be the quotient of $T \times H_X$
by the kernel of ineffectivity of its action on
$\rq{X}$. 
Then $G$ acts effectively on $\rq{X}$ and
Proposition~\ref{gencoxprops}~(ii) tells us 
that $\rq{X}$ is $G$-factorial.

Consider the $T$-invariant prime divisors
$E_1, \ldots, E_m \subseteq X$ 
supported in $X \setminus X_0$.
Their inverse images $\rq{E}_k := p^{-1}(E_k)$
are $G$-prime divisors and, 
since $\rq{X}$ is $G$-factorial, we have 
$\rq{E}_k = \div(f_k)$ with some $G$-homogeneous
$f_i \in \mathcal{R}(X)$.
According to Proposition~\ref{liftaction}~(ii),
the $\rq{E}_k$ are precisely the $G$-prime
divisors supported in $\rq{X} \setminus \rq{X}_0$.
Moreover, consider the $T$-invariant prime divisors 
$D_1, \ldots, D_n \subseteq X$  
along which $T$ acts with a  
finite generic isotropy group of order $l_j > 1$
and their ($G$-prime) inverse images 
$\rq{D}_j := p^{-1}(D_j)$.
As before, we see that 
$\rq{D}_j = \div(g_j)$ holds with some 
$G$-homogeneous $g_j \in \mathcal{R}(X)$
and the generic 
isotropy group of the $G$-action on  
$\rq{D}_j$ 
has order $l_j$.
Note that 
none of the $\rq{D}_j$ equals one of
the $\rq{E}_k$. 
Moreover, we may view the functions $f_k$ and $g_j$
as the canonical sections of the divisors 
$E_k$ and $D_j$.

Let $G_k \subseteq G$ denote the 
generic isotropy group of 
$\rq{E}_k \subseteq \rq{X}_k$.
The action of $G_0 := G_1 \cdots G_m$
on $\rq{X}_0 = p^{-1}(X_0)$
admits a geometric quotient 
$\lambda_0 \colon \rq{X}_0 \to \rq{Y}_0$.
The factor group $H := G/G_0$  
acts with at most finite isotropy groups
on $\rq{Y}_0$ and, 
by Proposition~\ref{prop:1isosplit}~(ii),
has generic isotropy group $H_j \subseteq H$
of order $l_j$ along the divisors 
$\rq{C}_j := \lambda_0(\rq{D}_j)$. 
Set $H_0 := H_1 \cdots H_n$ and 
let $\kappa \colon \rq{Y}_0 \to \rq{Z}_0$
denote the quotient for the action of 
$H_0$ on $\rq{Y}_0$.
The induced action of $F := H/H_0$ on 
$\rq{Z}_0$ admits again a geometric quotient 
$\rq{Z}_0 \to Z_0$ and the whole situation 
fits into the following 
commutative diagram.
\begin{equation}
\label{diag:x2x0modT}
\xymatrix{
{\rq{X}}
\ar@{}[r]|\supseteq
\ar[d]_{/ H_X}
& 
{\rq{X}_0}
\ar[r]^{/ G_0}
\ar[d]^{/ H_X}
& 
{\rq{Y}_0}
\ar[r]^{/ H_0}
&
{\rq{Z}_0}
\ar[d]^{/ F}
\\
X
\ar@{}[r]|\supseteq
& 
X_0
\ar[rr]_{/T}
& &
Z_0
}
\end{equation}
Replacing $Z_0$ and $X_0$ as well as $\rq{Z}_0$, $\rq{Y}_0$ 
and $\rq{X}_0$ with suitable 
big open subsets, we achieve that the group $F$ acts freely 
on $\rq{Z}_0$. 

We show that $\rq{Z}_0 \to Z_0$ is a universal
torsor for $Z_0$. 
According to Proposition~\ref{bigfreequot} this means 
to verify that $\rq{Z}_0$ is an $F$-factorial
quasiaffine variety with only constant invertible 
functions.
Proposition~\ref{prop:1isosplit} provides the 
corresponding properties for the $H$-variety
$\rq{Y}_0$.
Moreover, by Lemma~\ref{lem:fixinv}, every
$g_j$ is $G_0$-invariant, hence $g_j$ descends 
to a function on $\rq{Y}_0$, where we have 
$\div(g_j) = \rq{C}_j$.
Thus, we can apply Proposition~\ref{lem:quaff2free}
to obtain the desired properties for $\rq{Z}_0$
and the action of $F$. 

The final task is to relate the Cox rings
$\mathcal{R}(X) = \Gamma(\rq{X}, \mathcal{O})$
and 
$\mathcal{R}(Z_0) = \Gamma(\rq{Z}_0, \mathcal{O})$
to each other.
Note that we have canonical inclusions of graded algebras
$$ 
\Gamma(\rq{X}, \mathcal{O})
\ \supseteq \
\Gamma(\rq{X}_0, \mathcal{O})^{G_0}
\ = \ 
\Gamma(\rq{Y}_0, \mathcal{O})
\ \supseteq \
\Gamma(\rq{Y}_0, \mathcal{O})^{H_0}
\ = \ 
\Gamma(\rq{Z}_0, \mathcal{O}),
$$
where the first one is due to Proposition~\ref{prop:1isosplit}.
This allows us in particular to view
$\Gamma(\rq{Z}_0, \mathcal{O})$ as a graded subalgebra
of $\Gamma(\rq{X}, \mathcal{O})$.
The assertion now follows from 
Proposition~\ref{prop:1isosplit}~(iv) 
and Proposition~\ref{lem:quaff2free}~(i).
\end{proof}

In the above proof, we realized a big open subset
of $X_0/T$ as a quotient of a quasiaffine 
variety with only constant invertible functions 
by a free action of a diagonalizable group,
see the diagram~\ref{diag:x2x0modT}.
According to Proposition~\ref{bigfreequot}, this
allowed us to define a Cox ring for $X_0/T$.
Moreover, we use this now to show that $X_0/T$ 
admits a separation.

\begin{proposition}
\label{prop:sepex}
Let $\mathcal{X}$ be a normal quasiaffine 
variety with a free action of 
a diagonalizable group $H$.
Suppose that every invertible function on
$\mathcal{X}$ is constant 
and that $\mathcal{X}$ is $H$-factorial.
Then $X := \mathcal{X}/H$ admits a separation.
\end{proposition}

\begin{proof}
We first treat the case of
a certain toric variety.
Consider the standard action of $\TT^r = (\KK^*)^r$ 
on $\KK^r$,
let $\mathcal{Z} \subseteq \KK^r$ be the union
of all orbits of the big torus $\TT^r \subseteq \KK^r$
of dimension at least $r-1$, and let 
$H \subseteq \TT^r$ be a closed subgroup acting 
freely on $\mathcal{Z}$.
The fan $\Sigma$ of $\mathcal{Z}$ has the 
extremal rays of the positive 
orthant $\QQ_{\ge 0}^r$ as its maximal cones
and $Z := \mathcal{Z}/H$ is the toric prevariety
obtained by gluing the orbit spaces 
$\mathcal{Z}_{\varrho}/H$ along their common
big torus $T/H$, 
where $\mathcal{Z}_{\varrho } \subseteq \mathcal{Z}$ 
denotes the affine toric chart corresponding to 
$\varrho \in \Sigma$.
The embedding $H \to \TT^r$ corresponds to a 
surjection $\ZZ^r \to K$ of the respective 
character groups.
Let $P \colon \ZZ^r \to N$ be a map 
having $\Hom(K,\ZZ)$ as its kernel.
Then we obtain a canonical separation $Z \to Z'$ 
onto a toric variety $Z'$, the fan of which 
lives in $N$ and consists of the cones
$P(\varrho)$, where $\varrho \in \Sigma$.

In the general case, choose a finitely 
generated graded
subalgebra $A \subseteq \Gamma(\mathcal{X},\mathcal{O})$
such that we obtain an open embedding 
$\mathcal{X} \subseteq \b{X}$, where $\b{X} := \Spec \, A$.
Properly enlarging $A$, we may assume that 
it admits a system $f_1, \ldots, f_r$ 
of homogeneous generators such that each 
$\div(f_i)$ is $H$-prime in $\rq{X}$. 
Consider the $H$-equivariant closed 
embedding $\b{X} \to \KK^r$ 
defined by $f_1, \ldots, f_r$ 
and let 
$\mathcal{Z} \subseteq \KK^r$ be 
as above.
By construction,  
$\mathcal{U} := \mathcal{Z} \cap \mathcal{X}$
is a big $H$-invariant open subset 
of $\mathcal{X}$,
and we obtain 
a commutative diagram
$$ 
\xymatrix{
{\mathcal{U}}
\ar[r]
\ar[d]_{/H}
&
{\mathcal{Z}}
\ar[d]^{/H}
\\
U 
\ar[r]
&
Z
}
$$
where the induced map $U \to Z$ of quotients 
is a locally closed embedding and $Z$ is a toric prevariety.
Again by construction, the intersection of the invariant 
prime divisors of $Z$ with $U$ are prime divisors on $U$.
Consequently, the restriction of $Z \to Z'$ defines the desired
separation $U \to U'$  
\end{proof}

\begin{proof}{Proof of Theorem~\ref{fingenchar2}}
We prove the statement more generally for the 
cases that 
$\Gamma(X, \mathcal{O}) = \KK$ holds
or that 
$\Gamma(X, \mathcal{O}^*) = \KK^*$ holds 
and $\Cl(X)$ is free.
By Proposition~\ref{prop:sepex}, 
the orbit space $X_0/T$ admits a 
separation $\pi \colon X_0/T \to Y$.
According to Remark~\ref{rem:goodsep}, 
we may assume that there are prime divisors
$C_0, \ldots, C_r$ on $Y$ 
such that each $\pi^{-1}(C_i)$
is a disjoint union of prime divisors 
$C_{ij} \subseteq X_0/T$,
where $1 \le j \le n_i$, 
the map $\pi$ is an isomorphism over 
$Y \setminus (C_0 \cup \ldots \cup C_r)$
and all the $D_j$ occur among the 
divisors $D_{ij} := q^{-1}(C_{ij})$.
Then, according to Proposition~\ref{prop:sepcox}, 
we have
\begin{eqnarray*}
\mathcal{R}(X_0 / T)
& \cong & 
\mathcal{R}(Y)[\t{T}_{ij}; \; 0 \le i \le r, \, 1 \le j \le n_i]
\ / \ 
\bangle{\t{T}_i - 1_{C_i}; \; 0 \le i \le r},
\end{eqnarray*}
where the variables $\t{T}_{ij}$ correspond to
the canonical sections $1_{C_{ij}}$ 
and we define $\t{T}_i := \t{T}_{i1} \cdots \t{T}_{in_i}$.
Let $l_{ij} \in \ZZ_{\ge 1}$ denote the order
of the generic isotropy group of $D_{ij} = q^{-1}(C_{ij})$.
Then, by Theorem~\ref{fingenchar}, we have
\begin{eqnarray*}
\mathcal{R}(X)
& \cong &
\mathcal{R}(X_0/T)[S_1, \ldots, S_m; \,  T_{ij}] 
\ / \ \bangle{T_{ij}^{l_{ij}} - 1_{C_{ij}}},
\end{eqnarray*}
where the variables $T_{ij}$ correspond to the 
canonical sections $1_{C_{ij}}$;
note that $1_{C_{ij}}$ and $1_{D_{ij}}$ are 
identified for $l_{ij} = 1$. 
Putting these two presentations of Cox ring
together, 
we arrive at the assertion.
\end{proof}

\begin{proof}{Proof of Theorem~\ref{complexity1}}
We prove the statement more generally for the 
case that $\Gamma(X, \mathcal{O}) = \KK$ holds.
Since the $T$-action on $X$ is of complexity 
one, the orbit space $X_0/T$ is of dimension one
and smooth.
Moreover, using the diagram~\ref{diag:x2x0modT}
and Proposition~\ref{prop:1isosplit}, 
we see that $X_0/T$ 
admits only constant global functions
and has a finitely generated divisor class group.
It follows that $X_0/T$ is isomorphic to 
$\PP_1(A,\mathfrak{n})$, with $A = (a_0, \ldots, a_r)$
and $\mathfrak{n} = (n_0, \ldots, n_r)$ defined as 
in Theorem~\ref{complexity1}.
By Proposition~\ref{P1AnCox},
the Cox ring $\PP_1(A,\mathfrak{n})$ is given by
\begin{eqnarray}
\label{eqn:cox1}
\mathcal{R}(\PP_1(A,\mathfrak{n}))
&  \cong &
\KK[\t{T}_{ij}] \ /  \ \bangle{\t{g}_i; \; 0 \le i \le r-2},
\end{eqnarray}
where the variables $\t{T}_{ij}$ correspond to 
the canonical sections of points 
$a_{ij}$ in $X_0/T \cong \PP_1(A,\mathfrak{n})$.
Their inverse images $D_{ij} = q^{-1}(a_{ij})$ 
under $q \colon X_0 \to X_0/T$ are prime divisors 
with generic isotropy group of order 
$l_{ij} \ge 1$; note that $l_{ij} = 1$ is allowed.
Applying Theorem~\ref{fingenchar} gives
\begin{eqnarray}
\label{eqn:cox2}
\mathcal{R}(X)
& \cong & 
\mathcal{R}(\PP_1(A,\mathfrak{n}))[S_1,\ldots, S_m,T_{ij}] 
\ / \ 
\bangle{T_{ij}^{l_{ij}} - \t{T}_{ij}},
\end{eqnarray}
where the variables $S_i$ correspond to the 
divisors of $X$ having a one-dimensional generic isotropy 
group, the variables $T_{ij}$ are the canonical sections
of the divisors $D_{ij}$, and the $\t{T}_{ij}$ are 
identified with their pullbacks under $X_0 \to X_0/T$;
note that the pullback $q^*(a_{ij})$ equals $l_{ij}D_{ij}$.
Now, putting the descriptions~(\ref{eqn:cox1}) 
and~(\ref{eqn:cox2}) together gives the assertion. 
\end{proof}

\begin{remark}
Note that for factorial affine varieties with a complexity 
one torus action, D.~Panyshev observed in~\cite[Remark~2.12]{Pan}
a presentation of the algebra of global functions
by generators and trinomial relations.
\end{remark}

\section{Cox ring via polyhedral divisors}
\label{sec:crviapdiv}

In this section, we combine Theorem~\ref{fingenchar2}
with the description of algebraic torus actions
in terms of polyhedral divisors presented 
in~\cite{MR2207875} and~\cite{divfans} 
and provide a combinatorial approach
to the Cox ring of an algebraic variety 
with torus action.
We begin with a brief reminder on  
the language of polyhedral divisors.

In the sequel, $N$ is a free finitely generated 
abelian group, and $M = \Hom(N,\ZZ)$ is its dual.
The associated rational vector spaces are denoted 
by $N_\QQ := N \otimes \QQ$ and $M_\QQ := M \otimes \QQ$.
Moreover, $\sigma \subseteq N_\QQ$ is a pointed convex polyhedral 
cone, and $\omega \subseteq M_\QQ$ is its dual cone.
The relative interior of $\sigma$ is denoted by $\sigma^\circ$,
and if $\tau$ is a face of $\sigma$, then we write 
$\tau \preceq \sigma$.

We consider convex polyhedra $\Delta \subseteq N_\QQ$ 
admitting a decomposition $\Delta = \Pi + \sigma$ 
with a (bounded) polytope $\Pi \subseteq N_\QQ$;
we refer to $\sigma$ as the {\em tail cone\/} of $\Delta$
and refer to $\Delta$ as a {\em $\sigma$-polyhedron}.
With respect to Minkowski addition, 
the set $\Pol^+_{\sigma}(N)$
of all $\sigma$-polyhedra is 
a monoid with neutral element $\sigma$.
We consider also the empty set as an 
element of $\Pol^+_{\sigma}(N)$
and set 
$ \Delta + \emptyset := \emptyset + \Delta := \emptyset$.

We are ready to enter the description
of affine varieties with an action of the torus 
$T = \Spec \, \KK[M]$.
Let $Y$ be a normal variety and fix a pointed convex 
polyhedral cone $\sigma \subseteq N_\QQ$.
A {\em polyhedral divisor\/} on $Y$
is a formal finite sum
\begin{eqnarray*}
\D 
& = & 
\sum_Z \Delta_Z \cdot Z,
\end{eqnarray*}
where $Z$ runs over the prime divisors of $Y$ 
and the coefficients $\Delta_Z$ belong to 
$\Pol^+_{\sigma}(N)$; finiteness of the sum 
means that only finitely many coefficients 
$\Delta_Z$ differ from the tail cone $\sigma$.

The {\em locus\/} of a polyhedral divisor $\D$ 
on $Y$ is the open subset $Y(\D) \subseteq Y$ 
obtained by removing all prime divisors 
$Z \subseteq Y$ with $\Delta_Z = \emptyset$.
For every $u \in \omega \cap M$ we have the 
evaluation
\begin{eqnarray*}
\D(u)
& := & 
\sum_Z \min_{v \in \Delta_Z} \bangle{u ,v} \mal Z,
\end{eqnarray*}
which is an ordinary rational divisor living 
on $Y(\D)$.
We call the polyhedral divisor $\D$ on $Y$ 
{\em proper\/} if its locus is semiprojective,
i.e., projective over some affine variety, 
and its evaluations $\D(u)$, where $u \in \omega \cap M$,
have the following properties
\begin{enumerate}
\item 
$\D(u)$ has a base point free multiple,
\item 
$\D(u)$ is big for $u\in \omega^{\circ} \cap M$.
\end{enumerate}

\begin{remark}
\label{rem:ppcrit}
Suppose that we have $Y = \PP_n$,
and consider a polyhedral divisor 
$\D = \sum \Delta_Z \mal Z$.
The {\em degree\/} of $\D$ 
is the polyhedron
$$ 
\deg(\D)
\ := \ 
\sum_Z \Delta_Z \cdot \deg(Z)
\ \in \ 
\Pol^+_{\sigma}(N).
$$
It provides a simple criterion for properness: 
if $\deg(\D)$ is a proper subset of the 
tail cone of $\D$, 
then $\D$ is a proper polyhedral divisor,
see~\cite[Ex.~2.12.]{MR2207875}.
\end{remark}

By construction, every polyhedral divisor 
$\D$ on a normal variety $Y$ defines 
a sheaf $\mathcal{A}(\D)$ of $M$-graded 
$\mathcal{O}_Y$-algebras and its ring
$A(\D)$ of global sections:
$$ 
\mathcal{A}(\D)
\ := \ 
\bigoplus_{u \in \omega \cap M} \mathcal{O}(\D(u)),
\qquad\qquad
A(\D) 
\ := \
\Gamma(Y(\D), \mathcal{A}(\D)).
$$
Now suppose that $\D$ is proper.
Then~\cite[Thm.~3.1]{MR2207875}
guarantees that $A(\D)$ is a normal affine algebra.
Thus, we obtain an affine variety 
$X(\D) := \Spec \, A(\D)$,
which comes with an effective action of
the torus $T = \Spec \, \KK[M]$.
By~\cite[Thm.~3.4]{MR2207875},
every normal affine variety with an effective torus 
action is isomorphic to some $X(\D)$.

\begin{example}
\label{ex:ppdiv1}
Set $N = \ZZ^2$, let $\sigma \subseteq N_\QQ$ be the cone 
generated by the vectors $(1,1)$ and $(0,1)$, and
consider the $\sigma$-polyhedra 
$\Delta_0$ and $\Delta_{\infty}$
given as follows:
$$
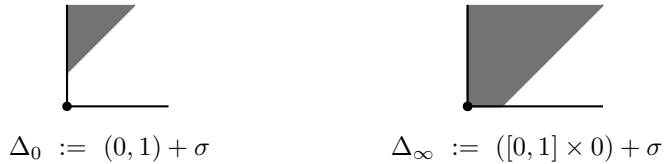
\begin{array}{ccc}
&&
\\
\psset{unit=0.45cm}
\begin{pspicture}(-.5,-.5)(3,2.5)%
  \psset{linecolor=black}%
  \psline[fillstyle=solid,linecolor=grayv,
          fillcolor=grayv]{-}(0,3)(0,1)(2,3)%
  \psline{-}(3,0)(0,0)(0,3)%
 \psset{framesep=0.5pt}
 \psdot(0,0)
 \color{black}
 \fontsize{8}{9}%
 \selectfont%
\end{pspicture}
&
{\hspace*{1.7cm}}
&
 \psset{unit=0.45cm}
 \begin{pspicture}(-.5,-.5)(4,2.5)%
 \psset{linecolor=black}%
 \psline[fillstyle=solid,linecolor=grayv,
         fillcolor=grayv]{-}(0,3)(0,0)(1,0)(4,3)%
 \psline{-}(4,0)(0,0)(0,3)%
 \psset{framesep=0.5pt}
 \psdot(0,0)
 \fontsize{8}{9}%
 \selectfont%
 \psdot(0,0)
 \color{black}
\end{pspicture}
\\
\Delta_0 
\ := \ 
(0,1) + \sigma
&
{\hspace*{1.7cm}}
&
\Delta_{\infty} 
\ := \ 
([0,1] \times 0) + \sigma
\end{array}
$$
Then we have a polyhedral 
divisor 
$\D := \Delta_0 \mal \{0\} + \Delta_\infty \mal \{\infty \}$
on $Y = \PP_1$.
Its degree $\deg(\D)$ and tail cone $\tail(\D)$ are 
given as
$$
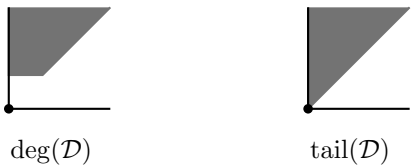
\begin{array}{ccc}
&&
\\
 \psset{unit=0.45cm}
 \begin{pspicture}(-.5,-.5)(3,2.8)%
   \psset{linecolor=black}%
 \psline[fillstyle=solid,linecolor=grayv,
         fillcolor=grayv]{-}(0,3)(0,1)(1,1)(3,3)
 \color{lgray}
 \psline{-}(3,0)(0,0)(0,3)%
 \psset{framesep=0.5pt}
 \psdot(0,0)
 
 \color{black}
 \fontsize{8}{9}%
 \selectfont%
\end{pspicture}
&
{\hspace*{1.7cm}}
&
 \psset{unit=0.45cm}
 \begin{pspicture}(-.5,-.5)(3,2.5)%
  \psset{linecolor=black}%
  \psline[fillstyle=solid,linecolor=grayv,
          fillcolor=grayv]{-}(0,3)(0,0)(3,3)%
 \psline{-}(3,0)(0,0)(0,3)%
 \psset{framesep=0.5pt}
 \psdot(0,0) 
 \color{black}
 \fontsize{8}{9}%
 \selectfont%
\end{pspicture}
\\
\deg(\D)
&
{\hspace*{1.7cm}}
&
\tail(\D)
\end{array}
$$
In particular, $\deg(\D)$ is a proper 
subset of $\tail(\D)$, and thus
Remark~\ref{rem:ppcrit} says that $\D$ 
is proper.
The associated $T$-variety is $\KK^3$ with the action 
\begin{eqnarray*} 
t \cdot z
& = & 
(t_1^{-1}t_2z_1,t_1z_2,t_2z_3).
\end{eqnarray*}
\end{example}

As in the case of toric varieties, 
general $T$-varieties are
obtained by gluing affine ones.
In the combinatorial picture, 
the gluing leads to the concept of 
a divisorial fan, which we recall
now.
As before, let $N$ be a finitely generated
free abelian group, fix a pointed convex
polyhedral cone $\sigma \in N_\QQ$, 
and let $Y$ be a normal variety.
Consider two polyhedral divisors
$$ 
\D
\ = \ 
\sum_Z \Delta_Z \cdot Z, 
\qquad \qquad 
\D'=\sum_Z \Delta'_Z \cdot Z
$$
both living on $Y$. 
The {\em intersection} of $\D$ and $\D'$ 
is the polyhedral divisor $\D \cap \D'$ on
$Y$ given by
\begin{eqnarray*}
\D \cap \D' 
& := & 
\sum_Z (\Delta'_Z \cap \Delta_Z) \cdot Z.
\end{eqnarray*}
Moreover, given a (not necessarily closed) point 
$y \in Y$, 
we define the {\em slice} of $\D$ at $y$
to be the polyhedron 
\begin{eqnarray*}
\D_y 
& := & 
\sum_{y \in Z} \Delta_Z.
\end{eqnarray*}
Note that the slice $\D_Y$ is the empty sum
and hence equals the tail cone of $\D$.
We say that $\D'$ is a {\em face\/} of $\D$
and write $\D' \preceq \D$
if $\D_y' \preceq \D_y$ holds for all $y \in Y$
and the $T$-equivariant morphism
$X(\D') \to X(\D)$ given by the inclusion
$A(\D') \supseteq A(\D)$ is an open embedding.

\begin{remark}
\label{rem:facecrit}
Suppose that in the above setting,
we have $Y = \PP_n$.
As a consequence of \cite[Lem. 6.7]{divfans} the relation $\D' \preceq \D$ holds
if and only if we have
$$
\D'_y \preceq \D_y \text{ for all } y \in Y,
\qquad\qquad
\deg(\D) \cap \tail(\D') \ = \ \deg(\D').
$$
\end{remark}

A {\em divisorial fan\/} is a finite set 
$\fan$ of polyhedral divisors such that for 
any two $\D,\D' \in \fan$ we have 
$\D \succeq \D' \cap \D \preceq \D'$.
For any $y \in Y$, we call the polyhedral complex 
$\fan_y$ defined by the slices $\D_y$ the
{\em slice\/} of $\fan$ at $y$. 
We say that the divisorial fan $\fan$ 
is {\em complete} if $Y$ is complete and each of its slices $\fan_y$ 
is a complete subdivision of $N_\QQ$.
The {\em locus} of $\fan$ is the 
open subset 
$$
Y(\fan) 
\ := \
\bigcup_{\D \in \fan} Y(\D) 
\ \subseteq \
Y.
$$

Given a divisorial fan $\fan$ consisting 
of proper polyhedral divisors, \cite[Thm.~5.3]{divfans} 
guarantees that we can 
equivariantly glue the affine $T$-varieties
$X(\D)$ along the open subsets $X(\D \cap \D')$,
where $\D, \D' \in \fan$, to a $T$-prevariety
$X(\fan)$.
If the divisorial fan $\fan$ is complete, 
then $X(\fan)$ is
a complete normal $T$-variety.
By~\cite[Thm.~5.6]{divfans}, every normal variety with 
torus action is isomorphic to some  $X(\fan)$.

\begin{example}
\label{sec:exmp-cotan}
Set $N := \ZZ^2$ and $Y := \PP_1$.
Consider the six polyhedral divisors
$\D^1, \ldots, \D^6$ with coefficients 
over the points $0, 1$ and $\infty$ as 
indicated below.
$$ 
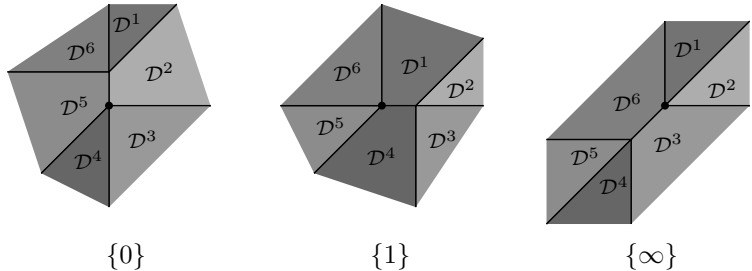
\begin{array}{ccccc}
&&&&
\\
 \psset{unit=0.45cm}
 \begin{pspicture}(-2.5,-3.7)(3.5,2.5)%
 \psset{linewidth=0.5pt}%
 \psline{-}(3,0)(0,0)(0,3)%
 \psline{-}(3,0)(0,0)(-2,-2)%
 \psline{-}(-2, -2)(0,0)(0,3)%
 \psline[fillstyle=solid,fillcolor=grayi]{-}(2,3)(0,1)(0,0)(3,0)%
 \psline[fillstyle=solid,fillcolor=grayii]{-}(3,0)(0,0)(0,-3)%
 \psline[fillstyle=solid,fillcolor=grayiii]{-}(-2, -2)(0,0)(0,1)(-3,1)%
 \psline[fillstyle=solid,fillcolor=grayiv]{-}(-3, 1)(0,1)(0,3)%
 \psline[fillstyle=solid,fillcolor=grayv]{-}(0,3)(0,1)(2,3)%
 \psline[fillstyle=solid,fillcolor=grayvi]{-}(0,-3)(0,0)(-2,-2)%
 \psset{framesep=0.5pt}
 \fontsize{6}{7}%
 \selectfont%
 \psdot(0,0)
 \color{black}
 \fontsize{8}{9}%
 \selectfont%
 \rput(0.55,2.45){$\D^1$}
 \rput(-0.6,-1.6){$\D^4$}
 \rput(1.5,1){$\D^2$}
 \rput(1,-1){$\D^3$}
 \rput(-1,0){$\D^5$}
 \rput(-0.8,1.6){$\D^6$}  
\end{pspicture}
&
&
 \psset{unit=0.45cm}
 \begin{pspicture}(-3,-3.7)(3.5,2.5)%
 \psset{linewidth=0.5pt}%
 \psline{-}(3,0)(0,0)(0,3)%
 \psline{-}(3,0)(0,0)(-2,-2)%
 \psline{-}(-2, -2)(0,0)(0,3)%
 \psline[fillstyle=solid,fillcolor=grayi]{-}(3,2)(1,0)(3,0)%
 \psline[fillstyle=solid,fillcolor=grayii]{-}(3,0)(1,0)(1,-3)%
 \psline[fillstyle=solid,fillcolor=grayiii]{-}(-2, -2)(0,0)(-3,0)%
 \psline[fillstyle=solid,fillcolor=grayiv]{-}(-3, 0)(0,0)(0,3)%
 \psline[fillstyle=solid,fillcolor=grayv]{-}(0,3)(0,0)(1,0)(3,2)%
 \psline[fillstyle=solid,fillcolor=grayvi]{-}(1,-3)(1,0)(0,0)(-2,-2)%
 \psset{framesep=0.5pt}
 \fontsize{6}{7}%
 \selectfont%
 \psdot(0,0)
 \fontsize{8}{9}%
 \selectfont%
 \color{black}
 \rput(1,1.2){$\D^1$}
 \rput(0,-1.5){$\D^4$}
 \rput(2.35,0.5){$\D^2$}
 \rput(1.8,-1){$\D^3$}
 \rput(-1.5,-0.6){$\D^5$}
 \rput(-1,1){$\D^6$}
\end{pspicture}
&
&
  \psset{unit=0.45cm}
  \begin{pspicture}(-3.3,-3.7)(2.5,2.5)%
    \psset{linewidth=0.5pt}%
    \psline{-}(2.5,0)(0,0)(0,2.5)%
    \psline{-}(2.5,0)(0,0)(-2.5,-2.5)%
    \psline{-}(-2.5, -2.5)(0,0)(0,2.5)%
    \psline[fillstyle=solid,fillcolor=grayi]{-}(2.5,2.5)(0,0)(0,0)(2.5,0)%
    \psline[fillstyle=solid,fillcolor=grayii]{-}(2.5,0)(0,0)(-1,-1)(-1,-3.5)%
    \psline[fillstyle=solid,fillcolor=grayiii]{-}(-3.5, -3.5)(-1,-1)(-3.5,-1)%
    \psline[fillstyle=solid,fillcolor=grayiv]{-}(-3.5, -1)(-1,-1)(0,0)(0,2.5)%
    \psline[fillstyle=solid,fillcolor=grayvi]{-}(-1,-3.5)(-1,-1)(-3.5,-3.5)%
    \psline[fillstyle=solid,fillcolor=grayv]{-}(0,2.5)(0,0)(2.5,2.5)%
    \psset{framesep=0.5pt}
    \fontsize{6}{7}%
    \selectfont%
    \psdot(0,0)
    \color{black}
    \fontsize{8}{9}%
    \selectfont
    \rput(0.7,1.7){$\D^1$}
    \rput(-1.5,-2.35){$\D^4$}
    \rput(1.7,0.5){$\D^2$}
    \rput(0.1,-1.1){$\D^3$}
    \rput(-2.35,-1.5){$\D^5$}
    \rput(-1.1,0.1){$\D^6$}
\end{pspicture}
\\
\{0\}
&
&
\{1\}
&
&
\{\infty\}
\end{array}
$$
The collections of degrees $\deg(\D^i)$ and tail cones 
$\tail(\D^i)$ of these polyhedral divisors are given as
$$ 
\begin{array}{ccc}
&&
\\
 \psset{unit=0.45cm}
 \begin{pspicture}(-2.5,-2.6)(3.5,2.8)%
  \psset{linewidth=.5pt}%
 \psline{-}(3,0)(0,0)(0,3)%
 \psline{-}(3,0)(0,0)(-2,-2)%
 \psline{-}(-2, -2)(0,0)(0,3)%
 \psline[fillstyle=solid,fillcolor=grayi]{-}(3,3)(1,1)(1,0)(3,0)%
 \psline[fillstyle=solid,fillcolor=grayii]{-}(3,0)(1,0)(0,-1)(0,-3)%
 \psline[fillstyle=solid,fillcolor=grayiii]{-}(-3, -3)(-1,-1)(-1,0)(-3,0)%
 \psline[fillstyle=solid,fillcolor=grayiv]{-}(-3, 0)(-1,0)(0,1)(0,3)%
 \psline[fillstyle=solid,fillcolor=grayv]{-}(0,3)(0,1)(1,1)(3,3)%
 \psline[fillstyle=solid,fillcolor=grayvi]{-}(0,-3)(0,-1)(-1,-1)(-3,-3)%
 \psset{framesep=0.5pt}
 \fontsize{6}{7}%
 \selectfont%
 \psdot(0,0)
 \color{black}
 \fontsize{8}{9}%
 \selectfont%
 \rput(1,2){$\D^1$}
 \rput(-0.6,-1.8){$\D^4$}
 \rput(2,1){$\D^2$}
 \rput(1,-1){$\D^3$}
 \rput(-1.5,-.5){$\D^5$}
 \rput(-0.8,1.3){$\D^6$}
\end{pspicture}
&
{\hspace*{1.7cm}}
&
 \psset{unit=0.45cm}
 \begin{pspicture}(-2.5,-2.6)(3.5,2.8)%
 \psset{linewidth=.5pt}%
 \psline{-}(3,0)(0,0)(0,3)%
 \psline{-}(3,0)(0,0)(-2,-2)%
 \psline{-}(-2, -2)(0,0)(0,3)%
 \psline[fillstyle=solid,fillcolor=grayi]{-}(3,3)(0,0)(3,0)%
 \psline[fillstyle=solid,fillcolor=grayii]{-}(3,0)(0,0)(0,-3)%
 \psline[fillstyle=solid,fillcolor=grayiii]{-}(-3, -3)(0,0)(-3,0)%
 \psline[fillstyle=solid,fillcolor=grayiv]{-}(-3, 0)(0,0)(0,3)%
 \psline[fillstyle=solid,fillcolor=grayv]{-}(0,3)(0,0)(3,3)%
 \psline[fillstyle=solid,fillcolor=grayvi]{-}(0,-3)(0,0)(-3,-3)%
 \psset{framesep=0.5pt}
 \fontsize{6}{7}%
 \selectfont%
 \psdot(0,0)
 \color{black}
 \fontsize{8}{9}%
 \selectfont%
 \rput(1,2){$\D^1$}
 \rput(-0.6,-1.8){$\D^4$}
 \rput(2,1){$\D^2$}
 \rput(1,-1){$\D^3$}
 \rput(-1.5,-.5){$\D^5$}
 \rput(-0.8,1.3){$\D^6$}
\end{pspicture}
\\
\\
\deg(\D^i)
&
{\hspace*{1.7cm}}
&
\tail(\D^i)
\end{array}
$$
Remarks~\ref{rem:ppcrit} and~\ref{rem:facecrit}
yield that 
$\D^1, \ldots, \D^6$ are proper and form a divisorial 
fan $\fan$.
The $T$-variety $X(\fan)$ is the projectivized 
cotangent bundle over $\PP_2$.
\end{example}

For the description of the Cox ring of 
the $T$-variety $X$ defined by a divisorial 
fan, we first of all need a description 
of the invariant prime divisors of $X$
and their generic isotropy groups.
For this, we introduce the following data. 

\begin{definition}
Consider a divisorial fan $\fan$ on a normal 
projective variety $Y$, and let $Z \subseteq Y$
be a prime divisor.
\begin{enumerate}
\item
The {\em index\/} of a vertex $v \in \fan_Z$ is 
the minimal positive integer $\mu(v)$ such that 
$\mu(v) \mal v \in N$ holds.
\item
We call a vertex $v \in \fan_Z$  {\em extremal}
if there is a $\D \in \fan$ with $v \in \D_Z$
such that $\mathcal{O}(\D(u))$ is big 
on $Z$ for any $u \in ((\D_Z -v)^{\vee})^{\circ}$.
The set of all extremal vertices $v \in \fan_Z$
is denoted by $\fan_Z^{\times}$
\item
We call a ray $\varrho \in \fan_Y$ {\em extremal}
if there is a $\D \in \fan$ with $\varrho \in \D_Y$
such that $\mathcal{O}(\D(u))$ is big 
on $Y$ for any $u \in (\varrho^{\perp}\cap\omega)^{\circ}$.
The set of all extremal rays $\varrho \in \fan_Y$
is denoted by $\fan_Y^{\times}$.
\item
We say that the prime divisor $Z$ is {\em irrelevant} if 
$\fan_Z^{\times}$ is empty,
and we denote by $Y^\circ \subseteq  Y(\fan)$ 
the open subset obtained by removing 
all irrelevant $Z$.
\end{enumerate}
\end{definition}

\begin{remark}
\label{rem:deg2extremal}
Let $\fan$ be a divisorial fan on $Y = \PP_n$ and 
$Z \subseteq \PP_n$ a prime divisor.
Then every vertex $v \in \fan_Z$ is extremal and 
a ray $\varrho \in \fan_Y$ is extremal
if and only if 
$\varrho \cap \deg(\D) = \emptyset$ holds
for some $\D \in \fan$ with $\varrho \in \D_Y$.
\end{remark}

As shown in~\cite{tidiv}, the extremal vertices 
of $\fan$ are in bijection with
the invariant prime divisors of $X = X(\fan)$
intersecting $X_0$ and the extremal rays 
correspond to those contained in 
$X \setminus X_0$; see also 
Propositions~\ref{sec:prop-primedivs-1} 
and~\ref{sec:prop-primedivs-2}.
We will denote by $D_v$ the divisors given 
by extremal vertices $v \in \fan_Z^{\times}$
and by $E_{\varrho}$ those
given by extremal rays $\varrho \in \fan_Y^{\times}$.
Then the divisor class group $\Cl(X)$
can be described as follows,
see~\cite[Cor.~3.17]{tidiv}.

\begin{proposition}
\label{prop:clgroup}
Let $\fan$ be a divisorial fan on $Y$
and set $X = X(\fan)$.
Then~$\Cl(X)$ is generated by the classes 
$[D_v]$, $v \in \fan_Z^{\times}$ and 
$[E_{\varrho}]$, $\varrho \in \fan_Y^{\times}$ 
and the image of a canonical homomorphism 
$\Cl(Y^{\circ}) \to \Cl(X)$.
The relations among these generators are
$$ 
\sum_{v\in \xvers_Z} \mu(v)D_{v}
\ = \
[Z],
\qquad
\qquad
\sum_{\varrho \in \xrays}  
\bangle{u, v_\varrho} E_\varrho  
+ 
\sum_Z \sum_{v \in \xvers_Z} \mu(v) \bangle{u,v} D_{v}
\ = \ 
0,
$$
where $Z$ runs through the prime divisors of $Y$,
$u$ runs through (a basis of) the lattice $M$ and 
$v_\varrho \in \varrho$ denotes the primitive 
lattice vector.
\end{proposition}

We are ready to compute the Cox ring 
of a $T$-variety $X = X(\fan)$ in terms 
of its defining divisorial fan $\fan$ and 
the projective variety $Y$ carrying $\fan$.
Let $Z_0, \ldots, Z_r \subseteq Y$
be the prime divisors having nontrivial
slices $\fan_{Z_0}, \ldots, \fan_{Z_r}$.

\begin{theorem}
\label{sec:thm-cox-divfan}
There is a $\Cl(X)$-graded inclusion of 
Cox rings $\mathcal{R}(Y^{\circ}) \to \mathcal{R}(X)$ 
and a $\Cl(X)$-graded isomorphism
\begin{eqnarray*}
\mathcal{R}(X)
& \cong &
\mathcal{R}(Y^\circ)[S_\varrho ,T_v; \; 
       \varrho \in \fan_Y^{\times}, \, 
       v \in \fan_{Z_i}^{\times}, \, 0 \le i \le r]
\ / \ 
\bangle{T^{\mu_i}- 1_{Z_i}; \; i=0,\ldots, r},
\end{eqnarray*}
where we set $T^{\mu_i} := \prod_{v \in \fan_{Z_i}^{\times}} T_v^{\mu(v)}$
and $1_{Z_i} \in \mathcal{R}(Y^{\circ})$ denotes the canonical
section of the prime divisor $Z_i \subseteq Y$.
The grading is given by $\deg T_v = [D_v]$ and 
$\deg S_\varrho = [E_\varrho]$.
\end{theorem}

As a direct consequence, we obtain the following 
description of the Cox ring of a $T$-variety 
of complexity one.

\begin{corollary}
\label{sec:cor-cox-cplx-one}
Let $\fan$  be a divisorial fan on $Y=\PP^1$ 
having non-trivial slices 
$\fan_{a_0}, \ldots, \fan_{a_r}$. 
Then the Cox ring of $X = X(\fan)$ is given by
$$
\KK[S_{\varrho}, T_{v}; \; \varrho \in \fan_Y^{\times}, 
                       \;  v \in \fan_{a_0}^{\times} 
                           \dot{\cup} \ldots \dot{\cup}\, 
                           \fan_{a_r}^{\times}]
\ / \ 
\biggl\langle 
\sum_{i=0}^r 
\beta_i T^{\mu_i}; \; 
\mathbf{\beta} \in \syz(\tilde{a}_0,\ldots,\tilde{a}_r) 
\biggr\rangle,
$$
where $\t{a}_i \in \KK^2$ represents $a_i \in \PP^1$,
we set $T^{\mu_i} := \prod_{v \in \fan_{a_i}^{\times}} T_v^{\mu(v)}$
and  $\syz(\t{a}_0,\ldots,\t{a}_r)$ is a basis 
for the space of linear relations among
$\t{a}_0,\ldots,\t{a}_r$.
\end{corollary}

Note that an appropriate choice of a basis for the space of 
linear relations among $\t{a}_0,\ldots,\t{a}_r \in \KK^2$
gives a trinomial representation of 
the Cox ring as in Theorem~\ref{complexity1}. 

\begin{example}
\label{sec:exmp-cot-1}
Consider once more the divisorial fan $\fan$ 
and its associated variety $X(\fan)$ of 
Example~\ref{sec:exmp-cotan}.
According to Remark~\ref{rem:deg2extremal},
there are no extremal rays and all 
six vertices
$$ 
v_1,v_2 \ \in \ \fan_{\{0\}},
\qquad
v_3,v_4 \ \in \ \fan_{\{1\}},
\qquad
v_5,v_6 \ \in \ \fan_{\{\infty\}}
$$
are extremal, where we have
$v_1=v_3=v_5=0 \in N$. 
Proposition~\ref{prop:clgroup} shows that
$\Cl(X(\fan))$ is freely generated by 
the classes of $D_{v_1}$ and $D_{v_2}$.
By Corollary~\ref{sec:cor-cox-cplx-one},
the Cox ring of $X(\fan)$ is 
\begin{eqnarray*}
\mathcal{R}(X(\fan))
& = &
\KK[T_1, \ldots, T_6] \ / \  \bangle{T_1T_2 + T_3T_4 + T_5T_6}
\end{eqnarray*}
with 
$\deg T_1 = \deg T_3 = \deg T_5 = [D_{v_1}]$
and 
$\deg T_2 = \deg T_4 = \deg T_6 = [D_{v_2}]$.
Note that this presentation of the Cox ring 
shows that $X(\fan)$ can be obtained as a 
$\KK^*$-quotient of the Grassmannian $G(2,4)$.
\end{example}

The rest of the section is devoted to proving
Theorem~\ref{sec:thm-cox-divfan},
which basically means to express the input 
data of Theorem~\ref{fingenchar2} in terms
of polyhedral divisors.
For this, we first have to recall further details 
of the construction of the $T$-variety 
$X(\fan)$ associated to a divisorial fan 
$\fan$ on a projective variety $Y$.
For every $\D \in \fan$, we have the sheaf
$\mathcal{A}(\D)$ of normal $M$-graded 
$\mathcal{O}_Y$-algebras.
Its relative spectrum
$\t{X}(\D) := \Spec_Y \mathcal{A}(\D)$
comes with a $T$-action
and we have canonical morphisms 
$$
\t{X}(\D) \ \to \ Y,
\qquad \qquad
\t{X}(\D) \ \to \ X(\D)
$$
defined by $\mathcal{O}_Y \subseteq \mathcal{A}(\D)$ 
and $A(\D)=\Gamma(\t{X}(\D),\Of)$.
The $T$-varieties glue together
$\t{X}(\D)$ along the open subsets $\t{X}(\D \cap \D')$
to a $T$-variety $\t{X}(\fan)$.
These gluings are compatible with the above
maps and one obtains a commutative diagram
$$ 
\xymatrix{
{\t{X}(\fan)}
\ar[rr]^{r}
\ar[dr]_{\t{\pi}}
&&
X(\fan)
\ar@{-->}[dl]^{\pi}
\\
& Y &
}
$$
where $r \colon \t{X}(\fan) \to X(\fan)$ is 
$T$-equivariant, birational and proper,
$\t{\pi} \colon \t{X}(\fan) \to Y$ is $T$-invariant
and the rational
map $\pi \colon X(\fan) \to Y$ is 
defined in codimension two.
Note that image of $\t{\pi}$ is given by
$$
\t{\pi}(\t{X}(\fan))
\ = \ 
\bigcup_{\D \in \fan} Y(\D) 
\ \subseteq \
Y.
$$

The next step is a precise description of the 
$T$-invariant prime divisors $X(\fan)$,
see also~\cite[Prop.~3.13]{tidiv}.
Consider an extremal vertex $v \in \fan_Z^{\times}$
of $\D \in \fan$, where $Z \subseteq Y(\D)$ 
is a prime divisor.
These data define a homogeneous ideal
\[
I_{v} 
\ := \
\bigoplus_{u \in \D_Y^\vee \cap M} 
\Gamma(Y, \Of(\D(u))) \cap \{f \in K(Y); \; \ord_Z(f) > -\bangle{u,v}\}
\ \subseteq \
\Gamma(X(\D), \mathcal{O}),
\]
which turns out to be a prime ideal of height one.
We define the corresponding prime divisor 
$D_{v} \subseteq X(\fan)$ to be the closure
of the zero set of $I_{v}$.

\begin{proposition}
\label{sec:prop-primedivs-1}
Set $X := X(\fan)$.
The assignment $v \to D_v$ induces a 
bijection between the  extremal vertices 
of $\fan$ and the invariant prime divisors of 
$X$ intersecting $X_0$. 
The extremal vertices of $\fan_Z$ correspond 
to the invariant prime divisors contained in 
$\b{\pi^{-1}(Z)}$ and 
the generic isotropy group of  
$D_v$ is cyclic of order~$\mu(v)$. 
\end{proposition}

\begin{proof}
We may restrict to the affine case. 
Consider a proper polyhedral divisor $\D$,
the corresponding sheaf of algebras 
$\mathcal{A}:=\mathcal{A}(\D)$
and its algebra of global sections 
$A:=A(\D)$.
First we calculate the ideal of 
$\b{\pi^{-1}(Z)}=r(\t{\pi}^{-1}(Z))$.
The inverse image ideal 
sheaf of $\Of(-Z)$ is given by
\begin{eqnarray*}
 \Of(-Z) \cdot \mathcal{A}
& = &
  \bigoplus_{u \in \D_Y^\vee \cap M} 
  \Of(\lfloor \D(u) \rfloor  - Z).
\end{eqnarray*}
The radical of the ideal
$\Gamma(Y,\Of(-Z) \cdot \mathcal{A}) \subseteq A$ 
is exactly the ideal we are looking for. 
It is given by
\begin{eqnarray*}
I_Z
& = & 
\bigoplus_{u \in \D_Y^\vee \cap M} 
\Gamma(Y, \Of(\D(u))) \cap \{f \in K(Y); 
\; \ord_Z(f) > - \min \bangle{u,\D_Z}\}
\\
& = & 
\bigcap_{v \in \D_Z} I_v.
\end{eqnarray*}
Note that we have $(I_v)_u = A_u$ 
if $\bangle{u,v} \notin \ZZ$. 
Denote by $\kappa \colon \t{Z}\rightarrow Z$
the normalization.
If $\psi \colon \tilde{Y} \rightarrow Y$ is a
desingularization, then we
have $A(\D) = A(\psi^*\D)$ and $I_v = \tilde{I}_{v}$, where
$ \tilde{I}_{v}$ is the ideal in $A(\psi^*\D)$ corresponding to the vertex $v$ in
$(\psi^*\D)_{f^{-1}_*Z}$. 
Hence, in the following we may assume that
$Y$ is smooth and thus every prime divisor is Cartier. 
Then, for the corresponding affine subschemes 
$V(I_v)$ we obtain the coordinate rings
\begin{eqnarray*}
A/I_v 
& = & 
\bigoplus_{u\in (\D_Z-v)^\vee \cap M_v} 
\Gamma(Y,\mathcal{A}_u)/\Gamma(Y,\mathcal{A}_u(-Z))
\\
& \subseteq &
\bigoplus_{u\in (\D_Z-v)^\vee \cap M_v} 
\Gamma (Y,\mathcal{A}_u/\mathcal{A}_u (-Z))
\\
& \cong &
\bigoplus_{u\in (\D_Z-v)^\vee \cap M_v} \Gamma (Z,\mathcal{A}_u|_Z)
\\
& \subseteq &
\bigoplus_{u\in (\D_Z-v)^\vee \cap M_v} 
\Gamma (\t{Z},\kappa^*(\mathcal{A}_u|_Z))
\\
& \cong  &
A(\D_v).
\end{eqnarray*}
Here, we write 
$\mathcal{A}_u(-Z) := \mathcal{A}_u \otimes \Of(-Z)$
as usual,
$M_v \subseteq M$ is the sublattice 
consisting of all $u \in M$ 
with $\bangle{u,v} \in \ZZ$ 
and $\D_v$ is a polyhedral divisor 
on $Z$ with tail cone 
$\sigma_v := \QQ_{\geq 0} \cdot (\D_Z-v)$ 
and lattice $M_v^* \supset N$ defined via the inclusion 
$\imath \colon Z \hookrightarrow Y$
as follows
\begin{eqnarray*}
\D_v 
& := & 
\sum_{W} (\D_W + \sigma_v) \cdot (\kappa \circ \imath)^*W.
\end{eqnarray*}
Note that  $(\kappa \circ \imath)^* Z$ is defined only 
up to linear equivalence as a divisor on 
$\t{Z}$ but every choice will give 
isomorphic algebras $A(\D_v)$, 
compare~\cite[Cor.~8.9]{MR2207875}. 
Our condition on the bigness of $\D(u)$ 
for $u \in \sigma_v^\vee$ implies that 
$\D_v$ is indeed proper for any extremal
$v$.
Hence, $X(\D_v)$ is irreducible 
and of dimension  $(\dim X -1)$
in this case. 
If $v$ is not extremal then $\D_v$ 
is the pullback of a proper polyhedral divisor on 
\begin{eqnarray*}
Y_v
& := & 
\text{Proj} \bigoplus_{k} \Gamma (Y,\Of(k \cdot \textstyle \sum_i \D_v(u_i)))
\end{eqnarray*}
for some $u_i \in (\sigma_v^\vee)^\circ$. 
But $Y_v$ is of smaller dimension than $Y$, 
since $\Of(\D_v(u_i)) \cong \Of(D(u_i))|_Z$ is not big. 
This implies, that $X(\D_v)$ is of dimension
$$
\dim Y_v + \dim T 
\ < \ 
\dim Y + \dim T 
\ = \ 
\dim X - 1.
$$

Since $\D(u)$ is semi-ample and $\mathcal{A}_u=\Of(\D(u))$ holds,
$\oplus_{k \geq 0}  H^1(\mathcal{A}_{k \cdot u}(-Z))$ 
is finitely generated as a module over the ring 
$\oplus_{k \geq 0} \Gamma(Y,\mathcal{A}_{k \cdot u})$. 
The long exact cohomology sequence
$$
0 
\rightarrow 
H^0(\mathcal{A}_{k \cdot u} (-Z)) 
\rightarrow  
H^0(\mathcal{A}_{k\cdot u}) 
\rightarrow 
H^0(\mathcal{A}_{k\cdot u}/(\mathcal{A}_{k\cdot u}(-Z)) 
\rightarrow  
H^1(\mathcal{A}_{k \cdot u}(-Z)) 
\rightarrow \ldots
$$
shows that 
$\oplus_{k \geq 0}  
\Gamma(Y,\mathcal{A}_{k\cdot u}/\mathcal{A}_{k \cdot u}  (-Z))$ 
is a finitely generated module over  the ring
$\oplus_{k \geq 0}  
\Gamma(Y,\mathcal{A}_{k \cdot u})/\Gamma(Y,\mathcal{A}_{k \cdot u}(-Z))$. 
The fact that 
$\Gamma(\t{Z},\kappa^*(\oplus_{u} \mathcal{A}_{k \cdot u}|_Z ))$ 
is finitely generated over 
$\Gamma(Z,\oplus_{u} \mathcal{A}_{k \cdot u}|_Z)$ 
follows from the properties of the normalization map. 
Thus, $A(\D_v)$ is finitely generated over $A/I_v$ 
and $D_v$ is the image of $X(\D_v)$ under a finite morphism $f$,
hence, $D_v$ is irreducible and of codimension one;
it is not hard to see that $f$ even is the normalization map.

The fact that all homogeneous functions of weights
$u \notin M_v$ vanish on $D_v$ implies, 
that $T$ does not act effectively on $D_v$ but with 
generic isotropy group $M/M_v \cong \ZZ/\mu(v)\ZZ$.
\end{proof}

Now take an extremal ray $\varrho \in \fan_Y^{\times}$ 
with $\varrho \in \D_Y$, where $\D \in \fan$.
Then define the associated invariant prime divisor 
$E_{\varrho}$ of $X(\fan)$ to be the closure 
of the zero set of $V(X(\D),I_{\varrho})$,
where $I_{\varrho}$ is the homogeneous prime ideal of 
height one given by
$$
I_{\varrho} 
\ := \ 
\bigoplus_{u \in \D_Y^\vee \setminus \varrho^\perp} \Gamma(Y,\Of(\D(u)))
\ \subseteq \
\Gamma(X(\D), \mathcal{O}).
$$

\begin{proposition}
\label{sec:prop-primedivs-2}
Set $X := X(\fan)$.
The assignment $\varrho \to E_{\varrho}$ induces a 
bijection between the set of extremal rays 
of $\fan$ and the invariant prime divisors of 
$X$ contained in $X \setminus X_0$.
\end{proposition}

\begin{proof}
For a polyhedral divisor $\D$,
the invariant prime divisors of $\t{X}$ 
contained in $\t{X}/\t{X}_0$ 
correspond to the prime ideal sheaves
given by not necessarily extremal rays
$\varrho \in \D_Y(1)$ as follows
\begin{eqnarray*}
\mathcal{I}_{\varrho} 
& := &
\bigoplus_{u \in \D_Y^\vee \setminus \varrho^\perp} \Of(\D(u)).
\end{eqnarray*}
This can be seen locally.
Consider an affine open subset $U \subset Y$  
such that $\D|_U$ is trivial. 
Then $\t{X}(\D|_U) \subseteq \t{X}(\D)$ 
is an open inclusion and we have 
\begin{eqnarray*}
A(\D|_U)
& = & 
\Gamma(U,\Of_Y)[D_Y^\vee \cap M].
\end{eqnarray*}
Now the claim follows from standard toric geometry,
since the considered prime divisors correspond 
to ideals $I \subset A(\D|_U)$ with $I \cap \Gamma(U,\Of_Y) = 0$.

The image under $r$ corresponds to the ideal 
$I_v=\Gamma(\mathcal{I}_{\varrho})$ 
and for the coordinate ring of the corresponding 
subvariety we obtain
$$
A(\D)/I_v
\ = \ 
\bigoplus_{u \in \varrho^\perp \cap \D_Y^\vee \cap M} \Gamma(Y,\D(u))
\ = \ 
A(\D_\varrho).
$$
Here, $\D_\varrho := \sum_{Z} p(\D_Z) \cdot Z$
is a polyhedral divisor on $Y$ with tail cone $p(D_Y)$ and lattice $p(N)$, 
where $p$ is the projection $N_\QQ \rightarrow N_\QQ/\QQ\cdot\varrho$. 

The fact that $\varrho$ is extremal ensures 
that $D_\varrho$ is proper, 
which in turn implies that $V(I_\varrho)$ 
has codimension one.
\end{proof}

\begin{proof}{Proof of Theorem~\ref{sec:thm-cox-divfan}}
We first construct big open subsets
$Y' \subseteq Y^{\circ}$ and $X' \subseteq X_0$.
The set $Y'$ is obtained by removing from 
$Y^{\circ}$ all the intersections 
$Z_i \cap Z_j$, where $0 \le i < j \le r$.
To define $X'$,
denote by $\t{E} \subseteq \t{X}$ the exceptional
locus of the contraction $r \colon \t{X} \to X$
and set 
$$ 
\t{X}'
\ := \ 
(\pi^{-1}(Y') \cap \t{X}_0) \setminus \t{E}
\ \subseteq \
\t{X}_0,
\qquad\qquad
X' 
\ := \ 
r(\t{X}')
\ \subseteq \
X_0.
$$

Then $\pi \colon \t{X}' \to Y'$ is surjective
and $r \colon \t{X}' \to X'$ is an isomorphism.
Moreover, the $T$-invariant map $\pi' := \pi \circ r^{-1}$ 
factors as 
$$ 
\xymatrix{
X' 
\ar[rr]^{\pi'}
\ar[dr]
&& 
Y'
\\
&
X'/T
\ar[ur]_{\varphi}
}
$$
Note that $\varphi$ is birational and injective.
Thus,  $\varphi$ is a local isomorphism and hence 
a separation for $X'/T$.

By Proposition~\ref{sec:prop-primedivs-1},
the prime divisors corresponding to the 
extremal vertices $v \in \fan_{Z_i}^{\times}$ 
are precisely the irreducible components of 
the inverse image $(\pi')^{-1}(Z_i)$, 
and their generic $T$-isotropy is of order 
$\mu(v)$.
Moreover, by Proposition~\ref{sec:prop-primedivs-2},
the prime divisors in $X \setminus X_0$
correspond to the extremal rays of $\fan$.
Now the assertion follows from
Theorem~\ref{fingenchar2}.
\end{proof}

\section{Applications and Examples}
\label{applexam}

We first note some algebraic properties
of the Cox ring of a variety with
complexity one torus action.
Recall the following concepts 
from~\cite[Def.~3.1]{Ha2}.
Let $K$ be a finitely generated abelian 
group and $R = \bigoplus_{w \in K} R_w$ 
any $K$-graded integral $\KK$-algebra with
$R^* = \KK^*$. 
\begin{enumerate}
\item
We say that a nonzero nonunit $f \in R$ 
is {\em $K$-prime} if it is homogeneous and 
$f \vert gh$ with homogeneous $g,h \in R$
always implies $f \vert g$ or $f \vert h$.
\item
We say that an ideal $\mathfrak{a} \subset R$
is {\em $K$-prime\/} if it is homogeneous and 
for any two homogeneous 
$f,g \in R$ with $fg \in \mathfrak{a}$ 
one has $f \in \mathfrak{a}$ or $g \in \mathfrak{a}$.
\item
We say that a homogeneous prime ideal 
$\mathfrak{a} \subset R$ has 
{\em $K$-height $d$\/} if $d$ is maximal 
admitting a chain 
$\mathfrak{a}_0 \subset \mathfrak{a}_1 \subset \ldots 
\subset \mathfrak{a}_d = \mathfrak{a}$
of $K$-prime ideals.
\item
We say that the ring $R$ is {\em factorially graded} if 
every  $K$-prime ideal of $K$-height 
one is principal.
\end{enumerate}

Now, let $X$ be a complete normal variety 
with finitely generated divisor class 
group and an algebraic torus action 
$T \times X \to X$ of complexity one. 
Then Theorem~\ref{complexity1} provides 
a presentation of the Cox ring of $X$ as 
\begin{eqnarray*}
\mathcal{R}(X)
& \cong & 
\KK[S_1, \ldots, S_m, T_{ij}; \; 0 \le i \le r, \, 1 \le j \le n_i] 
\ / \ 
\bangle{g_i; \; 0 \le i \le r-2},
\end{eqnarray*}
where the variables $S_j$ and $T_{ij}$ are homogeneous 
with respect to the $\Cl(X)$-grading
and the relations $g_i$ are $\Cl(X)$-homogeneous 
trinomials all having the same degree.

\begin{proposition}
\label{cplxone2complint}
The Cox ring $\mathcal{R}(X)$ is factorially 
$\Cl(X)$-graded.
In the presentation of Theorem~\ref{complexity1},
the generators $S_k$ and $T_{ij}$ define pairwise 
nonassociated $\Cl(X)$-prime elements 
and $\mathcal{R}(X)$ is 
a complete intersection.
\end{proposition}

\begin{proof}
The fact that $\mathcal{R}(X)$ is factorially 
$\Cl(X)$-graded holds for any complete variety
with a finitely generated Cox ring, 
use for example~\cite[Prop.~3.2]{Ha2}.
Moreover, the variables $S_k$ and $T_{ij}$ define 
pairwise nonassociated $\Cl(X)$-prime elements,
because their divisors are pairwise different 
$H_X$-prime divisors, where $H_X = \Spec \, \KK[\Cl(X)]$,
use again~\cite[Prop.~3.2]{Ha2}.

We show that $\mathcal{R}(X)$ is a complete 
intersection. 
This means to verify that $\mathcal{R}(X)$ is 
of dimension $m + n_0 + \ldots + n_r - (r-1)$.
Consider the torsor $\rq{X} \to X$, and recall
from the proof of Theorem~\ref{fingenchar},
diagram~\ref{diag:x2x0modT},
that we have quotients
$$ 
\rq{Y}_0 \ = \ \rq{X}_0 / G_0,
\qquad
\rq{Z}_0 \ = \ \rq{Y}_0 / H_0,
\qquad
X_0 / T \ \cong \ \rq{Z}_0 / F
$$
where $G_0$ is an $m$-dimensional torus acting 
freely,
$H_0$ is a finite group
and $F$ is a diagonalizable group acting freely
and having the rank of $\Cl(X_0/T)$ as its 
dimension.
In our situation, $X_0/T \cong \PP_1(A, \mathfrak{n})$
is of dimension one and, 
by Proposition~\ref{prop:sepcox} has a divisor 
class group of 
rank $n_0 + \ldots + n_r - r$.
Thus, the dimension of $\mathcal{R}(X)$
equals 
$$
\dim(\rq{X}_0)
\ = \ 
m + \dim(\rq{Y}_0)  
\ = \ 
m + \dim(\rq{Z}_0) 
\ = \ 
m + n_0 + \ldots + n_r - r + 1.
$$
\end{proof}

We come to geometric applications
of this observation.
Note that each complete normal 
variety $X$ with finitely generated 
divisor class group 
and a complexity one torus action 
is rational, because 
$X_0/T \cong \PP_1(A, \mathfrak{n})$ 
is so.
Thus, the varieties $X$ in question are 
precisely the 
complete normal rational ones 
with a torus action of complexity one.
If we impose additionally the 
condition that any two
points of $X$ admit a common affine 
neighbourhood, which holds e.g.
for projective $X$,
then 
Proposition~\ref{cplxone2complint}
and~\cite[Thm.~4.19]{Ha2}
ensure that $X$ arises from a 
``bunched ring'' $(R, \mathfrak{F}, \Phi)$,
see~\cite[Def.~3.3, Constr.~3.4]{Ha2},
where we may take  
$R = \mathcal{R}(X)$ and 
$\mathfrak{F} = (S_k,T_{ij})$.
This allows us to apply the results
provided in~\cite{Ha2}.

\begin{corollary}
\label{cor:embedd}
Let $X$ be a complete normal rational variety with an 
effective algebraic torus action $T \times X \to X$
of complexity one and suppose that any two
points of $X$ admit a common affine neighbourhood.
Then there exists a closed embedding 
$\imath \colon X \to X'$ into a toric 
variety $X'$ with big torus $T' \subseteq X'$ such 
that 
\begin{enumerate}
\item
$\imath \colon X \to X'$ is equivariant
with respect to a $T$-action on $X'$ given 
by a monomorphism $T \to T'$,
\item
the image $\imath(X) \subseteq X'$ intersects
$T'$ and is a complete intersection of 
$T$-invariant hypersurfaces of $X'$,
\item
for every $T'$-invariant prime divisor $D' \subseteq X'$,
the inverse image $\imath^{-1}(D') \subseteq X$ 
is a prime divisor,
\item
$\imath \colon X \to X'$ defines a pullback isomorphism 
$\imath^* \colon \Cl(X') \to \Cl(X)$ on the level of 
divisor class groups.
\end{enumerate}
\end{corollary}

\begin{proof}
Apply the construction of a toric embedding 
given in~\cite[Constr.~3.13 and Prop.~3.14]{Ha2}
to the defining bunched ring $(R, \mathfrak{F}, \Phi)$
of $X$, where $R = \mathcal{R}(X)$ and 
$\mathfrak{F} = (S_k,T_{ij})$,
and use the fact that the $S_k$ as well as
the $T_{ij}$ are homogeneous with respect to a 
lifting of the $T$-action to the torsor.
\end{proof}

Recall from the introduction that
$E_k \subseteq X$ are the prime divisors
supported in $X \setminus X_0$,
that $D_{ij} \subseteq X$ are prime 
divisors intersecting $X_0$ and lying 
over a point $a_{i} \in X_0/T$ and 
$l_{ij}$ is the order of the generic 
isotropy group of $T$ along $D_{ij}$.

\goodbreak 

\begin{corollary}
\label{cor:canondiv}
Let $X$ be a complete normal rational variety with an 
effective algebraic torus action $T \times X \to X$
of complexity one.
\begin{enumerate}
\item
The cone of divisor classes without fixed components
is given by
$$
\qquad \qquad
\bigcap_{1 \le k \le m}
\cone([E_s], [D_{ij}]; \; s \ne k) 
\ \cap \ 
\bigcap_{\genfrac{}{}{0pt}{}{0 \le i \le r}{1 \le j \le n_i}}
\cone([E_k], [D_{st}]; \; (s,t) \ne (i,j)).
$$
\item 
For any $0 \le i \le r$, one obtains a canonical 
divisor for $X$ by 
$$
\max(0,r-1) \cdot \sum_{j=0}^{n_i} l_{ij} D_{ij} 
\ - \ 
\sum_{k=1}^m E_k
\ -  \ 
\sum_{i,j} D_{ij}.
$$
\end{enumerate}
\end{corollary}

\begin{proof}
By~\cite[Thm.~2.3]{K3}, there is a 
small birational transformation $X \to X'$ 
with a projective $X'$.
As $X$ and $X'$ share the same Cox ring, we 
may assume that $X$ is projective.
The assertions then follow 
from~\cite[Prop.~4.1 and Prop.~4.15]{Ha2}.
\end{proof}

Note that~\cite[Thm.~3.19]{tidiv} provides 
an equivalent description of the canonical 
divisor in terms of the defining divisorial
fan.

The first non-trivial examples of torus 
actions of complexity one are $\KK^*$-surfaces.
Let us look at their Cox rings.
Orlik and Wagreich
associate in~\cite{OrWa} to any 
smooth complete $\KK^*$-surface
$X$ without elliptic fixed points a graph of the 
following shape:
$$
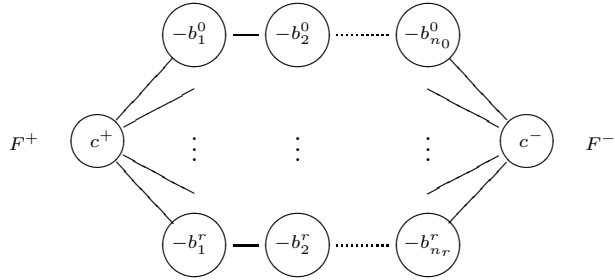
\entrymodifiers={++[o][F-]}
\def\objectstyle{\scriptstyle}
\xymatrix@-1.2pc{
*{}         
& *{}
& {-b^0_{1_{ \ }}} \ar@{-}[r]  
& {-b^0_{2_{ \ }}} \ar@{.}[rr] 
& *{}
& {-b^0_{n_0}}  
& *{}
& *{} 
\\
*{}
& *{}
& *{}
& *{}
& *{}
& *{}
& *{}
& *{}
\\
*{F^+}
&
{ \  c^+_{ \ }} 
\ar@{-}[uur] \ar@{-}[ur] \ar@{-}[ddr] \ar@{-}[dr] 
& *{\vdots}   
&*{\vdots}  
& *{}
& *{\vdots} 
& { \ c^-_{ \ }} 
\ar@{-}[ul] \ar@{-}[uul] \ar@{-}[ddl] \ar@{-}[dl] 
& *{F^-}
\\
*{}
& *{}
& *{}
& *{}
& *{}
& *{}
& *{}
& *{} 
\\
*{}
& *{}
& {-b_{1_{ \ }}^r} \ar@{-}[r]  
& {-b_{2_{ \ }}^r} \ar@{.}[rr] 
& *{}
& {-b_{n_r}^r} 
& *{}
& *{} 
\\
  }
$$
The vertices of this graph represent certain invariant 
curves.
The two (smooth) fixed point curves of $X$ occur as 
$F^+$ and $F^-$ in the graph.
The other vertices represent the invariant
irreducible contractible curves 
$D_{ij} \subseteq X$ different from $F^+$ and $F^-$.
The label $-b_{j}^i$ is the self intersection number
of $D_{ij}$, and two of the $D_{ij}$ are joined by an 
edge if and only if they have a common (fixed) point.
Every $D_{ij}$ is the closure of a non-trivial
$\KK^*$-orbit.

We show how to read off the Cox ring 
from the Orlik-Wagreich graph.
Suppose that $X$ is rational.
Then $F^+$ is rational as well
and hence is a $\PP_1$.
Define $l_{ij}$ to be the numerator of 
the canceled continued fraction
$$
b^i_1-\cfrac{1}{b^i_2-\cfrac{1}{\cdots -\cfrac{1}{b^i_{j-1}}}}
$$
Moreover, let $a_i$ be the point in $F^+ \cap D_{i1}$
and write $a_i = [b_i,c_i]$ with 
$b_i,c_i \in \KK$. 
Then, for every $0 \le i \le r$,
set 
$k = j+1 = i+2$ 
and define a trinomial
in $\KK[T_{ij}; \; 0 \le i \le r, \; 1 \le j \le n_i]$
as follows
$$
g_i 
\ := \
(c_kb_j - c_jb_k)f_i
\ + \ 
(c_ib_k - c_kb_i)f_j
\ + \ 
(c_jb_i - c_ib_j)f_k,
\quad
\text{where }
f_s 
\  := \
T_{s1}^{l_{s1}} \cdots T_{sn_s}^{l_{sn_s}}.
$$

\goodbreak

\begin{theorem}
\label{OrWa2Cox}
Let $X$ be a smooth complete rational $\KK^*$-surface
without elliptic fixed points.
Then the assignments $S^{\pm} \mapsto 1_{F^{\pm}}$ and 
$T_{ij} \mapsto 1_{D_{ij}}$ define an isomorphism
\begin{eqnarray*}
\mathcal{R}(X)
&  \cong &
\KK[S^+,S^-,T_{ij}; \;  \; 0 \le i \le r, \; 1 \le j \le n_i] 
\ / \ 
\bangle{g_i; \; 0 \le i \le r-2}
\end{eqnarray*}
of $\Cl(X)$-graded rings, where the $\Cl(X)$-grading on the 
right hand side is defined by associating to $S^{\pm}$ the 
class of $F^{\pm}$ and to $T_{ij}$ the class of $D_{ij}$.
\end{theorem}

\begin{proof}
The open set $X_0 \subseteq X$ is obtained by 
removing $F^+$, $F^-$ and the isolated 
fixed points.
By~\cite[Sec.~3.5]{OrWa}, the number $l_{ij}$ 
is the order of the isotropy group of 
the nontrivial $\KK^*$-orbit in $D_{ij}$.
Moreover,  we have a canonical morphism 
$\pi \colon X_0/\KK^* \to F^+$,
with exceptional fibers 
$\pi^{-1}(a_i) = \{a_{i1}, \ldots, a_{in_i}\}$,
where $a_{ij}$ represents the non-trivial 
$\KK^*$-orbit of $D_{ij}$. 
Thus, the assertion follows from Theorem~\ref{complexity1}.
\end{proof}

For (possibly singular) $\KK^*$-surfaces 
$X$ with elliptic fixed points,
the Cox ring can be computed as follows.
Suitably resolving gives a 
$\KK^*$-surface $\t{X}$, 
called canonical resolution, where the elliptic
fixed points are replaced with fixed point curves. 
Having computed the Cox ring 
$\mathcal{R}(\t{X})$ as above, 
we easily obtain the Cox ring 
$\mathcal{R}(X)$.
According to Theorem~\ref{complexity1}, we need 
the divisors of the type $E_k$ and $D_{ij}$ in $X$ and 
the orders $l_{ij}$ of the generic isotropy groups
of the $D_{ij}$.
Each of these divisors is the image of 
a non-exceptional divisor of the same 
type in $\t{X}$; to see this for the $D_{ij}$,
note that $X_0$ is the open subset of $\t{X}_0$
obtained by removing the exceptional locus 
of $\t{X} \to X$ and thus $X_0/\KK^*$
is an open subset of $\t{X}_0/\KK^*$.
Moreover, by equivariance, 
the orders $l_{ij}$ in $X$ are the same 
as in $\t{X}$. 
Consequently, the Cox ring $\mathcal{R}(X)$
is obtained 
from $\mathcal{R}(\t{X})$ by removing
those generators that correspond to the exceptional 
curves arising from the resolution.

As the intersection graphs of their 
resolutions are known, see~\cite{AlNi}, 
the methods just outlined provide 
Cox rings of (possibly singular)
Gorenstein del Pezzo $\KK^*$-surfaces~$X$;
note that Derenthal computed in~\cite{Der}
the Cox rings of the minimal resolutions~$\t{X}$ 
without assuming existence 
of a $\KK^*$-action for the cases 
that $X$ is of degree 
at least $3$ and $\mathcal{R}(\t{X})$ is defined 
by a single relation.
Moreover, the divisorial fans 
of Gorenstein del Pezzo $\KK^*$-surfaces $X$
are provided in~\cite{tfano},
which allows us to use as well 
the approach via 
polyhedral divisors.

\begin{example}
\label{ex:gordelpezzo}
We consider the family $X_\lambda$ of 
Gorenstein Del Pezzo $\KK^*$-surfaces 
over $\KK \setminus \{0,1\}$ of degree 
one and singularity type $\mathrm{2D_4}$. 
The canonical resolution~$\t{X}_\lambda$ of $X_\lambda$ 
is obtained by minimally resolving the two 
singularities 
and, by~\cite[Thm.~8.3]{AlNi}, its Orlik-Wagreich graph is 
given as 
$$
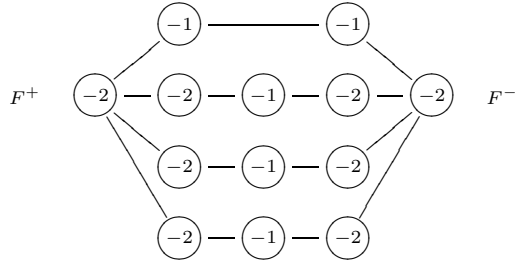
\entrymodifiers={++[o][F-]}
\def\objectstyle{\scriptstyle}
\xymatrix@-1.2pc{
*{}
&*{} 
& {-1} \ar@{-}[rr]  
&*{} 
& {-1} 
&*{} 
& *{} 
\\
*{F^+}
&{-2} \ar@{-}[r] \ar@{-}[ur] \ar@{-}[dr]\ar@{-}[ddr]  
& {-2} \ar@{-}[r] 
& {-1} \ar@{-}[r] 
& {-2} 
& {-2} \ar@{-}[ul] \ar@{-}[l] \ar@{-}[dl] \ar@{-}[ddl]
&*{F^-} 
\\
*{}         
&*{} 
& {-2} \ar@{-}[r]  
&{-1} \ar@{-}[r]
& {-2} 
&*{} 
& *{} 
\\
 *{}         
&*{} 
& {-2} \ar@{-}[r]  
&{-1} \ar@{-}[r]
& {-2} 
&*{} 
& *{} 
 }
$$
We have four points $a_0, \ldots, a_3$, where 
$\{a_i\} = D_{i1} \cap F^+$. 
Note that the positions of these four points on $F^+\cong \PP_1$ may 
vary and the parameter $\lambda$ is the cross ratio
of $a_0,a_1,a_2,a_3$. 
The Cox rings of $\t{X}_\lambda$ and $X_\lambda$ 
are given by 
\begin{eqnarray*}
\mathcal{R}(\t{X}_\lambda)
& = & 
\KK[S_1,S_2,T_{01},\ldots,T_{33}]
\ \biggl/  \ 
\biggl\langle
\begin{array}{l}
\scriptstyle
T_{01}T_{02}+T_{11}T_{12}^2T_{13}+T_{21}T_{22}^2T_{23},
\\ 
\scriptstyle
\lambda T_{11}T_{12}^2T_{13}+T_{21}T_{22}^2T_{23}+T_{31}T_{32}^2T_{33}
\end{array}
\biggr\rangle,
\\
\mathcal{R}(X_\lambda)
& = & 
\KK[T_1,\ldots,T_5]
\ \biggl/  \ 
\biggl\langle
\begin{array}{l}
\scriptstyle
T_1T_2 + T_3^2 + T_4^2,
\\ 
\scriptstyle
\lambda T_3^2 + T_4^2 + T_5^2
\end{array}
\biggr\rangle.
\end{eqnarray*}
Now let us look at $X_{\lambda}$ via its 
divisorial fan $\fan_{\lambda}$.
According to~\cite[Thm.~4.8]{tfano},
the divisorial fan $\fan_{\lambda}$
lives on $Y = \PP_1$. 
Its non-trivial slices lie over the 
points $a_0, \ldots, a_3 \in Y$ and 
are given in $N = \ZZ$ as follows:
$$
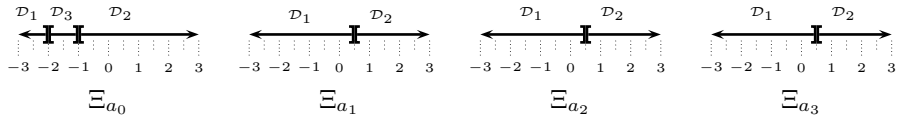
\begin{array}{cccc}
 \psset{unit=0.40cm}
 \begin{pspicture}(-3.4,-1)(3.4,1)%
   \fontsize{6}{7}%
   \mygrid
   \psset{linewidth=1pt}%
   \rput(-2.7,0.7){$\D_1$}
   \rput(-1.55,0.7){$\D_3$}
   \rput(.4,0.7){$\D_2$}
   \psline{<-]}(3,0)(-1,0)
   \psline{[-]}(-1,0)(-2,0)
   \psline{[->}(-2,0)(-3,0)
 \end{pspicture}
&
 \psset{unit=0.40cm}
 \begin{pspicture}(-3.4,-1)(3.4,1)%
  \fontsize{6}{7}%
  \mygrid
  \psset{linewidth=1pt}%
  \rput(-1.3,0.6){$\D_1$}
  \rput(1.4,0.6){$\D_2$}
  \psline{<-]}(-3,0)(0.5,0)
  \psline{[->}(0.5,0)(3,0)
 \end{pspicture}
&
 \psset{unit=0.40cm}
 \begin{pspicture}(-3.4,-1)(3.4,1)%
   \fontsize{6}{7}%
   \mygrid
   \psset{linewidth=1pt}%
   \rput(-1.3,0.7){$\D_1$}
   \rput(1.4,0.7){$\D_2$}
   \psline{<-]}(-3,0)(0.5,0)
   \psline{[->}(0.5,0)(3,0)
 \end{pspicture}
&
 \psset{unit=0.40cm}
 \begin{pspicture}(-3.4,-1)(3.4,1)%
  \fontsize{6}{7}%
  \mygrid
  \psset{linewidth=1pt}%
  \rput(-1.3,0.7){$\D_1$}
  \rput(1.4,0.7){$\D_2$}
  \psline{<-]}(-3,0)(0.5,0)
  \psline{[->}(0.5,0)(3,0)
 \end{pspicture}
\\[1ex]
\fan_{a_0}
&
\fan_{a_1}
&
\fan_{a_2}
&
\fan_{a_3}
\end{array}
$$
We compute the divisor class group
$\Cl(X_\lambda)$.
According to Remark~\ref{rem:deg2extremal}, we have 
two extremal vertices $v_1,v_2$ in $\fan_{a_0}$ and 
one extremal vertex $v_{i+2}$ in $\fan_{a_i}$ for $i=1,2,3$.
Let $D_i$ be the prime divisor associated to $v_i$ 
for $i = 1, \ldots, 5$ and denote by $D_0$
the positive generator of $\Cl(Y) = \ZZ$.
Then Proposition~\ref{prop:clgroup} tells us 
that the divisor class group $\Cl(X_\lambda)$ 
is $\ZZ D_0 \oplus \ldots \oplus  \ZZ D_5$
modulo the relations defined by the rows of the 
matrix
\begin{eqnarray*}
A 
& = & 
\left(
{
\begin{array}{rrrrrr}%
 -1& 1& 1& 0& 0& 0\\%
 -1& 0& 0& 2& 0& 0\\%
 -1& 0& 0 & 0& 2& 0\\%
 -1& 0& 0& 0& 0& 2\\%
 0& -2& -1& 1& 1& 1\\%
\end{array}
}
\right)
\end{eqnarray*}
The Smith Normal Form 
$S =  U \cdot A \cdot V$ 
with unimodular
transformation matrices $U$ and $V$
is given as
\begin{eqnarray*}
S
& = &  
\left(
{
\begin{array}{rrrrrr}
  1&0&0&0&0&0\\%
  0&1&0&0&0&0\\%
  0&0&1&0&0&0\\%
  0&0&0&2&0&0\\%
  0&0&0&0&2&0\\%
\end{array}
}
\right)
\end{eqnarray*}
In particular, we conclude 
$\Cl(X_{\lambda}) \cong \ZZ \oplus \ZZ/2\ZZ \oplus \ZZ/2\ZZ$.
Moreover, computing $V^{-1}$ we see that the class 
of $D_4$ generates the free part 
and the classes of $D_3 - D_5$ and $D_4 - D_5$ 
generate the cyclic parts.
Consulting Theorem~\ref{sec:thm-cox-divfan} gives the 
Cox ring 
\begin{eqnarray*}
\mathcal{R}(X_\lambda)
& = & 
\KK[T_1, \ldots, T_5]
\ / \ 
\bangle{
T_1T_2 + T_3^2 + T_4^2, \, \lambda T_3^2 + T_4^2 + T_5^2
}
\end{eqnarray*}
whith the grading  
$$ 
\deg(T_1)
\ = \
\deg(T_2) 
\ = \
(1,\b{1},\b{0}),
\qquad
\deg(T_3)
\ = \ 
(1,\b{1},\b{1}),
$$
$$
\deg(T_4)
\ = \ 
(1,\b{0},\b{0}),
\qquad
\deg(T_5)
\ = \ 
(1,\b{0},\b{1}).
$$
\end{example}

Proceeding as in this example, we 
are able to compute the 
Cox rings of all Gorenstein del Pezzo 
$\KK^*$-surfaces and their minimal 
resolutions. 
Here comes the result for the cases of 
Picard number one and two.

\begin{theorem}
\label{thm:crgordelps}
Let $X$ be a Gorenstein del Pezzo surface of 
Picard number at most two admitting a 
nontrivial $\KK^*$-action.
The following table provides the Cox rings
of~$X$ and its minimal resolution $\t{X}$ 
ordered by the degree $\deg(X)$ and 
the singularity type ${\rm S}(X)$.
{\small
\begin{center}
\begin{longtable}[htbp]{lll}
\multicolumn{3}{c}{\bf $\deg(X) = 1$} 
\\[1ex]
\midrule
${\rm S}(X)$ 
&  
$\mathcal{R}(X)$ 
& 
$\mathcal{R}(\t{X})$
\\
\cmidrule{1-3}
$\mathrm{2D_4}$ 
& 
${\scriptstyle \KK[ T_1,\ldots,T_5] \big/ 
\left\langle
      \begin{smallmatrix}
        T_1T_2+T_3^2+T_4^2,\\
        \lambda T_3^2+T_4^2+T_5^2 
      \end{smallmatrix}
\right\rangle}$ 
& 
${\scriptstyle \KK[S_1,S_2,T_1,\ldots,T_{11}] \big/
\left\langle
      \begin{smallmatrix}
        T_1T_2+T_6T_7T_3^2+T_8T_9T_4^2,\\
        \lambda T_6T_7T_3^2+T_8T_9T_4^2+T_{10}T_{11}T_5^2 
      \end{smallmatrix}
    \right\rangle}$ 
\\
\cmidrule{1-3}
$\mathrm{E_6A_2}$ 
& 
${\scriptstyle \KK[T_1,\ldots,T_4]/ 
 \langle T_1^2T_2+T_3^3+T_4^3\rangle}$ 
&  
    ${\scriptstyle \KK[ S, T_1,\ldots,T_{11}]/\langle 
    T_5T_1^2T_2+T_6T_7T_8^2T_3^3+T_9T_{10}T_{11}^2T_4^3\rangle}$
\\
\cmidrule{1-3}
$\mathrm{E_7A_1}$ 
& ${\scriptstyle \KK[ T_1,\ldots,T_4]/\langle 
    T_1^3T_2+T_3^4+T_4^2\rangle}$ 
&  
    ${\scriptstyle\KK[S, T_1,\ldots,T_{11}]/\langle 
         T_5T_6^2T_1^3T_2+T_7T_8^2T_9^3T_{10}^3T_3^4+T_{11}T_4^2\rangle}$
\\
\cmidrule{1-3}
$\mathrm{E_8}$ 
& ${\scriptstyle \KK[T_1,\ldots,T_4]/\langle 
    T_1^5T_2+T_3^3+T_4^2\rangle}$ 
&  
    ${\scriptstyle \KK[S, T_1,\ldots,T_{11}]/\langle 
        T_5T_6^2T_7^3T_8^4T_1^5T_2+T_9T_{10}^2T_3^3+T_{11}T_4^2\rangle}$
\\
\bottomrule 
\\[2ex]
\multicolumn{3}{c}{\bf $\deg(X) = 2$} 
\\[1ex]
\midrule
${\rm S}(X)$ 
&  
$\mathcal{R}(X)$ 
& 
$\mathcal{R}(\t{X})$
\\
\cmidrule{1-3}
$\mathrm{2A_3A_1}$ 
& ${\scriptstyle \KK[T_1,\ldots,T_4]/\langle 
    T_1T_2+T_3^2+T_4^2\rangle}$ 
&  
    ${\scriptstyle \KK[S_1,S_2, T_1,\ldots,T_9]/\langle 
    T_5T_1T_2+T_6T_7T_3^2+T_8T_9T_4^2\rangle}$
\\
\cmidrule{1-3}
$\mathrm{A_5A_2}$ 
& ${\scriptstyle \KK[T_1,\ldots,T_4]/\langle 
    T_1T_2+T_3^3+T_4^3\rangle}$ 
&  
    ${\scriptstyle \KK[S, T_1,\ldots,T_{10}]/\langle 
    T_1T_2+T_5T_6T_7^2T_3^3+T_8T_9T_{10}^2T_4^3\rangle}$
\\
\cmidrule{1-3}    
$\mathrm{D_43A_1}$ 
& 
${\scriptstyle \KK[S_1,T_1,T_2,T_3]/\langle 
    T_1^2+T_2^2+T_3^2\rangle}$ 
&  
    ${\scriptstyle \KK[S_1,S_2, T_1,\ldots,T_9]/\langle 
    T_4T_5T_1^2+T_6T_4T_2^2+T_8T_9T_3^2\rangle}$
\\
\cmidrule{1-3}
$\mathrm{D_6A_1}$ 
& ${\scriptstyle \KK[T_1,\ldots,T_4]/\langle 
    T_1^2T_2+T_3^4+T_4^2\rangle}$ 
&  
    ${\scriptstyle \KK[S, T_1,\ldots,T_{10}] /\langle 
    T_5T_1^2T_2+T_6T_7T_8^2T_9^3T_3^4+T_{10}T_4^2\rangle}$
\\
\cmidrule{1-3}
$\mathrm{E_7}$ 
& 
${\scriptstyle \KK[T_1,\ldots,T_4]/\langle 
    T_1^4T_2+T_3^3+T_4^2\rangle}$ 
&  
${\scriptstyle \KK[S, T_1,\ldots,T_{10}]/\langle 
    T_5T_6^2T_7^3T_1^4T_2+T_8T_9^2T_3^3+T_{10}T_4^2\rangle}$
\\
\cmidrule{1-3}    
$\mathrm{2A_3}$ 
& ${\scriptstyle \KK[T_1,\ldots,T_6] \big/ \left\langle
      \begin{smallmatrix}
        T_1T_2+T_3T_4+T_5^2,\\
        \lambda T_3T_4+T_5^2+T_6^2 
      \end{smallmatrix}
    \right\rangle}$ 
& $ {\scriptstyle \KK[S_1,S_2,T_1,\ldots,T_{10}] \big/ \left\langle
      \begin{smallmatrix}
        T_1T_2+T_3T_4+T_7T_8T_5^2,\\
        \lambda T_3T_4+T_7T_8T_5^2+T_{9}T_{10}T_6^2 
      \end{smallmatrix}
    \right\rangle}$ 
\\
\cmidrule{1-3}    
$\mathrm{D_5A_1}$ 
& ${\scriptstyle \KK[T_1,\ldots,T_5]/\langle 
    T_1T_2^2+T_3T_4^2+T_5^3\rangle}$ 
&  
    ${\scriptstyle \KK[S, T_1,\ldots,T_{10}]/\langle 
    T_6T_1T_2^2+T_7T_3T_4^2+T_8T_9T_{10}^2T_5^3\rangle}$
\\
\cmidrule{1-3}    
$\mathrm{E_6}$ 
& ${\scriptstyle \KK[T_1,\ldots,T_5]/\langle 
    T_1T_2^3+T_3T_4^3+T_5^2\rangle}$ 
&  
    ${\scriptstyle \KK[S, T_1,\ldots,T_{10}]/\langle 
    T_6T_7^2T_1T_2^3+T_8T_9^2T_3T_4^3+T_{10}T_5^2\rangle}$
\\
\bottomrule 
\\[2ex]
\multicolumn{3}{c}{\bf $\deg(X) = 3$} 
\\[1ex]
\midrule
${\rm S}(X)$ 
&  
$\mathcal{R}(X)$ 
& 
$\mathcal{R}(\t{X})$
\\
\cmidrule{1-3}
$\mathrm{A_5A_1}$ 
& 
${\scriptstyle \KK[T_1,\ldots,T_4]/\langle 
    T_1T_2+T_3^4+T_4^2\rangle}$ 
&  
    ${\scriptstyle \KK[S, T_1,\ldots,T_9]/\langle 
    T_1T_2+T_5T_6T_7^2T_8^3T_3^4+T_9T_4^2\rangle}$
\\
\cmidrule{1-3}
$\mathrm{E_6}$ 
& 
${\scriptstyle \KK[T_1,\ldots,T_4]/\langle 
    T_1^3T_2+T_3^3+T_4^2\rangle}$ 
&  
    ${\scriptstyle \KK[S, T_1,\ldots,T_9]/\langle 
    T_5T_6^2T_1^3T_2+T_7T_8^2T_3^3+T_9T_4^2\rangle}$
\\
 \cmidrule{1-3}   
$\mathrm{2A_2A_1}$ 
& 
${\scriptstyle \KK[T_1,\ldots,T_5]/\langle 
    T_1T_2+T_3T_4+T_5^2\rangle}$ 
&  
    ${\scriptstyle \KK[S_1,S_2, T_1,\ldots,T_8]/\langle 
    T_6T_1T_2+T_3T_4+T_7T_8T_5^2\rangle}$
\\
\cmidrule{1-3}    
$\mathrm{A_32A_1}$ 
& 
${\scriptstyle \KK[S_1,T_1,\ldots,T_4]/\langle 
    T_1T_2+T_3^2+T_4^2\rangle}$ 
&  
    ${\scriptstyle \KK[S_1,S_2, T_1,\ldots,T_8]/\langle 
    T_1T_2+T_5T_6T_3^2+T_7T_8T_4^2\rangle}$
\\
\cmidrule{1-3}    
$\mathrm{A_4A_1}$ 
& 
${\scriptstyle \KK[T_1,\ldots,T_5]/\langle 
    T_1T_2^2+T_3T_4+T_5^3\rangle}$ 
&  
    ${\scriptstyle \KK[S, T_1,\ldots,T_9]/\langle 
    T_6T_1T_2^2+T_3T_4+T_7T_8T_9^2T_5^3\rangle}$
\\
\cmidrule{1-3}    
$\mathrm{D_5}$ 
& 
${\scriptstyle \KK[T_1,\ldots,T_5]/\langle 
    T_1T_2^3+T_3T_4^2+T_5^2\rangle}$ 
&  
    ${\scriptstyle \KK[S, T_1,\ldots,T_9]/\langle 
    T_6T_7^2T_1T_2^3+T_8T_3T_4^2+T_9T_5^2\rangle}$
\\
\bottomrule
\\[2ex] 
\multicolumn{3}{c}{\bf $\deg(X) = 4$} 
\\[1ex]
\midrule
${\rm S}(X)$ 
&  
$\mathcal{R}(X)$ 
& 
$\mathcal{R}(\t{X})$
\\
\cmidrule{1-3}
$\mathrm{D_5}$ 
& 
${\scriptstyle \KK[T_1,\ldots,T_4]/\langle 
    T_1^2T_2+T_3^3+T_4^2\rangle}$ 
&  
    ${\scriptstyle \KK[S, T_1,\ldots,T_8]/\langle 
    T_5T_1^2T_2+T_6T_7^2T_3^3+T_8T_4^2\rangle}$
\\
 \cmidrule{1-3}   
$\mathrm{A_3A_1}$ 
& 
${\scriptstyle \KK[T_1,\ldots,T_5]/\langle 
    T_1T_2+T_3T_4+T_5^3\rangle}$ 
&  
    ${\scriptstyle \KK[S, T_1,\ldots,T_8]/\langle 
    T_1T_2+T_3T_4+T_6T_7T_8^2T_5^3\rangle}$
\\
 \cmidrule{1-3}   
$\mathrm{A_4}$ 
& 
${\scriptstyle \KK[T_1,\ldots,T_5]/\langle 
    T_1T_2^3+T_3T_4+T_5^2\rangle}$ 
&  
    ${\scriptstyle \KK[S, T_1,\ldots,T_8]/\langle 
    T_6T_7^2T_1T_2^3+T_3T_4+T_8T_5^2\rangle}$
\\
\cmidrule{1-3}
$\mathrm{D_4}$ 
& 
${\scriptstyle \KK[T_1,\ldots,T_5]/\langle 
    T_1T_2^2+T_3T_4^2+T_5^2\rangle}$ 
&  
    ${\scriptstyle \KK[S, T_1,\ldots,T_8]/\langle 
    T_6T_1T_2^2+T_7T_3T_4^2+T_8T_5^2\rangle}$
\\
   \bottomrule
\\[2ex] 
\multicolumn{3}{c}{\bf $\deg(X) = 5$} 
\\[1ex]
\midrule
${\rm S}(X)$ 
&  
$\mathcal{R}(X)$ 
& 
$\mathcal{R}(\t{X})$
\\
 \cmidrule{1-3}   
$\mathrm{A_3}$ 
& 
${\scriptstyle\KK[ T_1,\ldots,T_5]/\langle 
   T_1T_2^2+T_3T_4+T_5^2\rangle}$ 
&  
    ${\scriptstyle \KK[S, T_1,\ldots,T_7]/\langle 
    T_6T_1T_2^2+T_3T_4+T_7T_5^2\rangle}$
\\
\cmidrule{1-3}
$\mathrm{A_4}$ 
& 
${\scriptstyle \KK[T_1,\ldots,T_4]/\langle 
    T_1T_2+T_3^3+T_4^2\rangle}$ 
&  
    ${\scriptstyle \KK[S, T_1,\ldots,T_{7}]/\langle 
    T_1T_2+T_5T_6^2T_3^3+T_7T_4^2\rangle}$
\\
\bottomrule
\\[2ex] 
\multicolumn{3}{c}{\bf $\deg(X) = 6$} 
\\[1ex]
\midrule
${\rm S}(X)$ 
&  
$\mathcal{R}(X)$ 
& 
$\mathcal{R}(\t{X})$
\\
\cmidrule{1-3}
$\mathrm{A_2}$ 
& 
${\scriptstyle \KK[T_1,\ldots,T_5]/\langle 
   T_1T_2+T_3T_4+T_5^2\rangle}$ 
&  
    $\KK[{\scriptstyle S, T_1,\ldots,T_6}]/\langle 
    {\scriptstyle T_1T_2+T_3T_4+T_6T_5^2}\rangle$
\\
\bottomrule 
\end{longtable}
\end{center}
}
\end{theorem}

\goodbreak

Finally, we consider equivariant vector bundles over 
a toric variety $X$ arising from a fan $\Sigma$
and ask for the Cox rings of their projectivizations.
We will use Klyachko's description~\cite{klyachko:eqbundles}
of equivariant reflexive sheaves over $X$; 
we will follow Perling's notation~\cite{perling}
in terms of families of complete {\em increasing\/} 
filtrations.

We recall the basic constructions.
Let $\mathcal{E}$ be an equivariant
reflexive sheaf of rank $r$ 
on $X$. 
Then $\mathcal{E}$ is trivial over the big
torus $T \subseteq X$.
Moreover, for every ray $\varrho \in \Sigma^{(1)}$,
the sheaf $\mathcal{E}$ splits over the 
affine chart $X_{\varrho} \subseteq X$ 
and hence is even trivial there.
This gives us identifications
$$
\Gamma(X_{\varrho}, \mathcal{E}) 
\ \subseteq \
\Gamma(T, \mathcal{E}) 
\ =  \
E \otimes \Gamma(T, \Of_X) 
\ =  \ 
E \otimes \KK[M]
$$
with an $r$-dimensional vector space $E$.
Fix generators 
$e_{\varrho,1} \otimes \chi^{u_{\varrho,1}},
\ldots, e_{\varrho,r} \otimes \chi^{u_{\varrho,r}}$
for every 
$\Gamma(X_{\varrho}, \mathcal{E})$.
Then $\mathcal{E}$ is 
determined by the family of complete
increasing filtrations $E^\varrho(i)$,
where $\varrho \in \Sigma^{(1)}$, 
of $E$ defined by 
\begin{eqnarray*}
E^\varrho(i)
& := &
\lin(e_{\varrho, j} ; \; \bangle{u_{\varrho,j}, v_\varrho} \le i),
\end{eqnarray*}
where $v_\varrho \in \varrho$ denotes 
the primitive lattice vector.
Conversely, given any family of complete
increasing filtrations $E^\varrho(i)$,
where $\varrho \in \Sigma^{(1)}$, 
of $E = \KK^r$, one obtains an equivariant 
reflexive sheaf $\mathcal{E}$ of rank $r$
over $X$ by defining its sections over 
the affine charts $X_{\sigma} \subseteq X$,
where $\sigma \in \Sigma$, to be 
$$
\Gamma(X_\sigma, \mathcal{E})
\ := \
\bigoplus_{u \in M} 
\left(
\bigcap_{\varrho \in \sigma^{(1)}} 
E^\varrho (\bangle{u, v_\varrho})
\right) 
\otimes 
\chi^u 
\ \subseteq \
E \otimes \KK[M].
$$

In our first result, we compute the Cox ring 
of the projectivization $\PP(\mathcal{E})$,
see~\cite[p.~162]{Ht}, of a locally free 
sheaf $\mathcal{E}$ of rank two
over a complete toric variety
$X$ arising from a fan $\Sigma$.
Let $E^\varrho(i)$, 
where $\varrho \in \Sigma^{(1)}$, 
be the family of filtrations
describing $\mathcal{E}$,
let $\mathcal{L}$ be the set of 
one-dimensional subspaces of $E$ 
occurring in these filtrations,
for every $L \in \mathcal{L}$ fix 
a generator $e_L$,
and denote by $\syz(\mathcal{L})$ 
the space of linear relations among 
the $e_L$.
Moreover, let $i_k^\varrho$ be the smallest integer 
such that $\dim E^\varrho(i_k^\varrho) > k$
and set $L^\varrho:=E^\varrho(i_0^\varrho)$. 

\goodbreak

\begin{theorem}
\label{sec:cor-bundle}
Let $\mathcal{E}$ be an equivariant 
locally free sheaf of rank two
over a complete toric variety~$X$ defined by a
fan $\Sigma$. 
Then the Cox ring of the projectivization 
$\PP(\mathcal{E})$ is given as
\begin{eqnarray*}
\mathcal{R}(\PP(\mathcal{E})) 
& = &
\KK[S_\varrho, T_L; \; \varrho \in \Sigma^{(1)}, \, L \in \mathcal{L}]
\ \Bigl/ \
\left\langle 
\sum_{L \in \mathcal{L}} \lambda_L S^L T_{L}; \; 
\mathbf{\lambda} \in \syz\left(\mathcal{L}\right) 
\right\rangle,
\\
& & 
\text{where}
\quad 
S^L 
\ := \ 
\prod_{\varrho,\; L^\varrho=L} S_\varrho^{i^\varrho_1-i^\varrho_0}.
\end{eqnarray*}
\end{theorem}

\begin{example}
Let $\mathcal{T}$ be the sheaf of sections 
of the tangent bundle of the projective plane 
$\PP_2$; then $\PP(\mathcal{T})$ is the 
projectivized cotangent bundle.  
As a toric variety, $\PP_2$ is given by the 
complete fan in $\QQ^2$ with the rays 
$$
\varrho_1 \ =\ \QQ_{\geq 0} \mal e_1,
\qquad
\varrho_2 \ = \ \QQ_{\geq 0}\mal e_2,
\qquad
\varrho_0 \ = \ \QQ_{\geq 0}\mal e_0,
$$
where $e_1,e_2 \in \QQ^2$ are the canonical 
basis vectors and we set $e_0 := -e_1-e_2$.
The  filtrations of the tangent
sheaf are given as
$$
E^\varrho(i)=
\begin{cases}
0, & \; i < -1,
\\
\KK \mal \varrho, & \; i = -1,
\\
E, & \; i > -1.
\end{cases}
$$
As generators for the one-dimensional subspaces 
we may choose 
$e_1, e_2, e_0 \in \KK^2$. 
The linear relations between them are spanned by 
$(1,1,1) \in \KK^3$. Hence, as in
Example~\ref{sec:exmp-cot-1}, we obtain 
\begin{eqnarray*}
\mathcal{R}(\PP(\mathcal{T}))
& = & 
\KK[S_{1},S_{2}, S_{3}, T_{1}, T_2, T_3] 
\ / \  
\bangle{T_1S_1 + T_2S_2 + T_3S_3}.
\end{eqnarray*}
\end{example}

More generally, we may calculate the Cox ring 
of the projectivized cotangent bundle on an arbitrary  
smooth complete toric variety $X$ arising from a fan $\Sigma$.
We distinguish two types of rays $\varrho \in \Sigma^{(1)}$:
those with $-\varrho \notin \Sigma^{(1)}$ and 
those with $-\varrho \in \Sigma^{(1)}$. 
Denote by $\mathcal{L}$ the set containing 
all rays of the first type 
and one representative for every pair 
of the second type. 
Moreover, let 
$\syz(\mathcal{L})$ denote the tuples
$\mathbf{\lambda} \in \KK^{\mathcal{L}}$ 
such that 
$\sum_{\varrho \in \mathcal{L}} \lambda_\varrho v_\varrho = 0$,
where $v_\varrho \in \varrho$ denote the primitive generator. 

\begin{theorem}
\label{sec:cor-cotan}
Let $X$ be a smooth complete toric variety arising 
from a fan $\Sigma$, and denote by $\mathcal{T}_X$
the sheaf of sections of the tangent bundle over $X$.
Then the Cox ring of the projectivization
$\PP(\mathcal{T}_X)$ is given by
\begin{eqnarray*}
\mathcal{R}(\PP(\mathcal{T}_X))
& = & 
\KK[S_\varrho, T_\tau; \; \varrho \in \Sigma^{(1)}, \, \tau \in \mathcal{L}]
\ \Bigl / 
\left\langle 
\sum_{\varrho \in \mathcal{L}} \lambda_\varrho S^\varrho T_\varrho; \; 
\mathbf{\lambda} \in \syz\left(\mathcal{L}\right)  
\right\rangle,
\\
& & 
\text{where} \quad 
S^\varrho \ := \
\begin{cases}
  S_\varrho S_{-\varrho} & -\varrho \in \Sigma^{(1)},
\\
  S_\varrho & else.
\end{cases}
\end{eqnarray*}
\end{theorem}

\begin{remark}
\label{sec:rem-toric-quot-sep}
Let $\Sigma$ be a fan in a lattice $N$
having rays $\varrho_1, \ldots, \varrho_s$ 
as its maximal cones,
$Z$ the associated 
toric variety and 
$T \subseteq Z$ 
the acting torus.
Given a primitive sublattice $L \subseteq N$, 
consider the action of the corresponding 
subtorus $H \subseteq T$ on~$Z$.
Let $P \colon N \to N' := N/L$ denote the 
projection.
The generic isotropy group 
$H_{\varrho} \subseteq H$ 
along the toric divisor $D_{\varrho} \subseteq Z$ 
corresponding to a ray $\varrho \in \Sigma$ 
is one-dimensional if $P(\varrho) = 0$ holds and 
finite otherwise; in the latter case it is 
given by
\begin{eqnarray*}
\Chi(H_{\varrho})
& = & 
(\lin(P(\varrho)) \cap N') 
\ / \ 
P(\lin(\varrho) \cap N).
\end{eqnarray*}

In particular, the set 
$Z_0 \subseteq Z$ is the toric 
subvariety corresponding to the subfan 
$\Sigma_0 \subseteq \Sigma$
obtained by removing all $\varrho$ with 
$P(\varrho) = \{0\}$.
For an affine chart $Z_{\varrho} \subseteq Z_0$,
the orbit space $Z_{\varrho} / H$ is the 
affine toric variety  $Z'_{P(\varrho)}$ 
corresponding to the ray $P(\varrho)$ in $N'$.
Gluing these $Z'_{P(\varrho)}$ along 
their common big torus $T/H$ gives the 
toric prevariety $Z_0 / H$.

There is a canonical separation 
$\pi \colon Z_0 / H \to Z'$, 
where $Z'$ is the toric variety 
defined by the fan $P(\Sigma_0)$ in $N'$
having $\{P(\varrho); \; \varrho \in \Sigma_0\}$
as its set of maximal cones.
Note that the inverse image $\pi^{-1}(D_{P(\varrho)})$
of the divisor $D_{P(\varrho)} \subseteq Z'$ 
corresponding to $P(\varrho) \in P(\Sigma_0)$
is the disjoint union of all divisors 
$D_{\tau} \subseteq Z_0$ with  
$P(\tau) = P(\varrho)$.
\end{remark}

\begin{proof}{Proof of Theorem~\ref{sec:cor-bundle} and Theorem~\ref{sec:cor-cotan}}
In order to use Theorem~\ref{fingenchar2},
we have to study the map $\pi\circ q$ 
obtained by composing the quotient 
$q \colon \PP(\mathcal{E})_0 \rightarrow \PP(\mathcal{E})_0/T$ 
with the separation 
$\pi \colon \PP(\mathcal{E})_0/T \rightarrow Y$. 
This done in three steps. 
First we cover $\PP(\mathcal{E})$ by affine 
toric charts and describe the quotient map 
on these charts using Remark~\ref{sec:rem-toric-quot-sep}. 
Then we collect the data for 
Theorem~\ref{fingenchar2} in every chart. 
In the last step we will see how this local 
data fit into the global picture.

\medskip\noindent 
{\em Step 1: the toric charts.} 
\enspace
We may assume that the maximal cones 
of $\Sigma$ are just the rays 
$\varrho \in \Sigma$.  
On an affine chart $X_\varrho$ any equivariant 
locally free sheaf $\mathcal{E}$ is actually free 
with homogeneous generators 
$s_{\varrho,0}, \ldots, s_{\varrho,r}$ 
of the form 
$s_{\varrho,i} = e_{\varrho,i} \otimes \chi^{u_{\varrho,i}}$. 
Choosing an appropriate order, we may achieve
$\bangle{u_{\varrho,k},v_\varrho}=i_k$.
Then 
$\PP(\mathcal{E}|_{X_\varrho})$ is given as
$$
\PP(\mathcal{E}|_{X_\varrho}) 
\ = \ 
\Proj_{X_\varrho}(S(\mathcal{E}|_{X_\varrho})) 
\ = \ 
\Proj\; \KK[\varrho^\vee \cap M][s_{\varrho,0}, \ldots, s_{\varrho,r}], 
$$ 
where $\deg(s_{\varrho,i}):=1$.
So, $\PP(\mathcal{E}|_{X_\varrho})$ is 
$X_\varrho \times \PP^r$ but endowed 
with a special $T^n$-action. 
This action can be extended to an 
$T^{n+r}$-action by assigning to $s_i$ 
the weight $(u_i,b_i) \in M \times M'$ 
for $i=0, \dots r$. 
Here, $M' \cong \ZZ^r$ and $b_1, \ldots, b_r$ 
is a basis and $b_0:=0$. 
As a consequence we can describe 
$\PP(\mathcal{E}|_{X_\varrho})$ 
as a toric variety and the 
endowed $T^n$-action by an 
inclusion $T^n \hookrightarrow T^{n+r}$,
which corresponds to the lattice 
inclusion $N \hookrightarrow N \times N'$, 
where $N':=(M')^*$. 

We describe the corresponding fan 
in $N \times N'$.
Denote by  $b^*_1, \ldots, b^*_r$ the dual 
basis in $N'$ and set $b^*_0 := -\sum b^*_i$.
Moreover, set
$$
\varrho_i  
\ := \
\tilde{\varrho} + 
\cone(b^*_0, \ldots , b^*_{i-1}, b^*_{i+1}
, \ldots, b^*_r),
\quad
\tilde{\varrho}
\ := \
\QQ_{\geq 0}\cdot (v_\varrho, - \sum_{k=0}^r i^\varrho_k b^*_k).
$$
Then $\varrho_0, \ldots, \varrho_r$ are the 
maximal cones of the fan we are looking for. 
Indeed, cover $\PP(\mathcal{E}|_{X_\varrho})$ 
by the affine charts 
$\Spec\; \KK[\varrho^\vee \cap M][\frac{s_0}{s_i}, \ldots, \frac{s_r}{s_i}]$. 
Then 
\begin{eqnarray*}
\KK[\varrho^\vee \cap M][\textstyle\frac{s_0}{s_i}, \ldots, \frac{s_r}{s_i}]  
& \to & 
\KK[\varrho_i^\vee \cap M \times M'],
\\
\frac{s_j}{s_i} 
& \mapsto & 
\chi^{(u_j-u_i,b_j-b_i)},
\\
\chi^{u} 
& \mapsto & 
\chi^{(u,0)}
\end{eqnarray*}
defines an equivariant isomorphism 
from $\PP(\mathcal{E}|_{X_\varrho})$ with 
the extended torus action 
onto the toric variety arising from the 
fan just defined, see also~\cite[pp. 58-59]{oda:cb}. 

\medskip
\noindent 
{\em Step 2: local quotient maps.} 
\enspace
In the setting of Remark~\ref{sec:rem-toric-quot-sep},
the map $P \colon N \times N' \rightarrow N'$ 
is the projection to the second factor 
and we deduce that the separation 
of $\PP(\mathcal{E}|_{X_\varrho})_0/T$ 
has the fan $\Sigma'$ consisting of the rays 
$\QQ_{\geq 0} \cdot b^*_0, \ldots, \QQ_{\geq 0} \cdot b^*_r, 
P(\tilde{\varrho})$ and the trivial cone.

If $\mathcal{E}$ is an equvariant locally free 
sheaf of rank two, then we obtain 
$\tilde{\varrho}=\QQ_{\geq 0} \cdot (v_\varrho, (i_1-i_0)b_0)$ 
and $\Sigma'$ is the unique fan $\Delta$ of $\PP^1$. 
If $\mathcal{E}$ is the sheaf of sections 
of the tangent bundle of $X$, then we have $E=N\otimes \KK$ 
and the filtrations
$$
E^\varrho(i)=
\begin{cases}
0, & \; i < -1,
\\
\KK \cdot \varrho, & \; i = -1,
\\
E, & \; i > -1.
\end{cases}
$$
Thus, by the chosen order of the $s_{\varrho,i}$,
we obtain 
$\tilde{\varrho}=\QQ_{\geq 0} \cdot (v_\varrho, b_0)$ 
and 
$\Sigma'=\{0,\QQ_{\geq 0} \cdot b^*_0, \ldots, \QQ_{\geq 0} \cdot b^*_r\}$. 
Hence, $\Sigma'$ is a subfan of a fan $\Delta$
with $X_\Delta \cong \PP^r$ 
and we have a rational toric map 
$p_\varrho:\PP(\mathcal{E}|_{X_\varrho}) \dashrightarrow \PP^r$,
which is defined on a big open subset.

Now we locally collect the data for 
Theorem~\ref{fingenchar2} using 
Remark~\ref{sec:rem-toric-quot-sep}. 
In the preimage $P^{-1}(\QQ_{\ge 0} \cdot b_0)$,
we find the rays $\tau:=\QQ_{\ge 0}\cdot (0, b_0)$ 
and $\tilde{\varrho}$. 
Hence, the prime divisor in $Y=X_{\fan'}$ 
corresponding to $\QQ_{\ge 0} \cdot b_0$ 
has two invariant prime divisors in its 
preimage under the map $\pi \circ q$.

The lattice elements of $\tau$ are mapped 
onto the lattice elements 
of $\QQ_{\geq 0} \cdot b_0$, 
hence $T$ acts effectively on the 
corresponding prime divisor. 
The lattice generated by $P(\tilde{\varrho}\cap (N\times N'))$ 
has index $i^\varrho_1-i^\varrho_0$ in $\ZZ \cdot b_0$. 
This implies that the corresponding prime divisor 
has a generic isotropy group of order $i_1-i_0$ 
(which is equal to $1$ in the case of $\mathcal{T}_X$).

\medskip
\noindent 
{\em Step 3: the global picture.} 
\enspace
We identify $\PP(E^*)$ with $X_\Delta$ via the 
isomorphism $\varphi_\varrho$ induced by the 
following homomorphism of the homogeneous coordinate rings
$$
S(E) \ \to \ \KK[\chi^{b_0}, \ldots, \chi^{b_r}],
\qquad  
e_i  \ \mapsto \ \chi^{b_i}.
$$
Note that the map $\varphi_\varrho \circ p_\varrho$ 
 (considered as a rational map $\PP(\mathcal{E}) \dashrightarrow \PP(E^*)$)
no longer depends on the choice of $\varrho$ and 
$s_{\varrho,0}, \ldots, s_{\varrho,r}$, 
because over $\PP(\mathcal{E}|_T)$ it is just
the projection 
$\PP(\mathcal{E}|_T)=\PP(E^*)\times T \rightarrow \PP(E^*)$ 
given by 
\begin{equation*}
  \begin{array}{rcccl}
    S(E) 
& \hookrightarrow &  
\KK[M][s_{\varrho,0}, \ldots, s_{\varrho,r}] 
& = & 
\Gamma(T, S(\mathcal{E}))
\\
e_i  
&  \mapsto &  
\chi^{-u_i} \cdot s_{\varrho,i} 
& = & 
e_i \otimes \chi^0
\end{array}.
\end{equation*}
Note that $\varphi_\varrho$ maps the prime divisor 
corresponding to the ray $\QQ_{\geq 0}\cdot b_0$ 
onto $E^\varrho(i_0)^\perp=e_{\varrho,0}^\perp \subset \PP(E^*)$.
Putting things together we obtain
\begin{itemize}
\item 
A ray $\varrho$ with $i^\varrho_0 = i^\varrho_1$ 
corresponding to a divisor outside of $\PP(\mathcal{E})_0$.
\item 
Since $(\pi \circ q)|_{\PP(\mathcal{E}|_T)}$ 
is the equivariant projection 
$\PP(\mathcal{E}|_T)=\PP(E^*)\times T \rightarrow \PP(E^*)$ 
the closure $\b{(\pi \circ q)^{-1}(Z) \cap \PP(\mathcal{E}|_T)}$ 
is always a prime divisor with effective $T$-action.
\item 
If $i^\varrho_0 < i^\varrho_1$ the ray $\varrho$ corresponds 
to an additional invariant prime divisor in the preimage 
of $E^\varrho(i_0)^\perp$ with generic isotropy group of order $i_1-i_0$.
\end{itemize}
Inspecting the filtrations for $\mathcal{T}_X$,
we see that  $i^\varrho_1 - i^\varrho_0=1$ 
for every ray $\varrho$ and 
$E^\varrho(i^\varrho_0) = E^\tau(i^\tau_0)$ 
if and only if $\tau = \pm\varrho$. 
Using Theorem~\ref{fingenchar2}, 
we obtain the desired results.
\end{proof}

\end{document}